\documentstyle{amsppt}
\magnification 1200
\vsize=8.9 true in

\define \ZZ {\Bbb Z}
\define \RR {\Bbb R}
\define \NN {\Bbb N}
\define \QQ {\Bbb Q}
\define \CC {\Bbb C}
\define \EE {\Bbb E}
\define \KK {\Bbb K}
\define \FF {\Bbb F}
\define \HH {\Bbb H}
\define \SS {\Bbb S}
\redefine \AA {\Bbb A}
\define \OO {\Bbb O}
\define \agcheck{\alpha^{\vee}}
\define \bgcheck{\beta^{\vee}}
\define \Wtilde{\wt W}
\define \Stilde{\wt S}
\define \Tbar{\bar T}
\define \Sgbar{\bar \Sigma}
\define \Hscr{\Cal H}
\define \What{\hat W}

\define \ac {\acute}
\define \fa {\forall}
\define \id {\operatorname{id}}

\define \sign{\operatorname{sign}}
\define \ord{\operatorname{ord}}
\define \Sym{\operatorname{Sym}}
\define \Ind{\operatorname{Ind}}
\define \End{\operatorname{End}}
\define \Hom{\operatorname{Hom}}
\define \Irr{\operatorname{Irr}}
\define \Lie{\operatorname{Lie}}
\define \gr{\operatorname{gr}}
\define \Mod{\text{\bf Mod}}
\define \sgn{\operatorname{sgn}}
\define \Tr{\operatorname{Tr}}


\define\arr{\rightarrow}
\define\arriso{@> {\thicksim\ } >>}
\define\asr{\mapsto}


\define \wt{\widetilde}
\define \w0M{{w_{0,M}}}


\newcount\parno
\newcount\secno
\newcount\subsecno
\global\parno=0
\global\secno=0
\global\subsecno=0

\define\newparagraph#1{
   \global\secno=0
   \global\advance\parno by 1
   \bigpagebreak\bigpagebreak\medskip
   \subheading\nofrills{\the\parno.~#1}
   \par\nobreak}

\define\secstart{
   \global\subsecno=0
   \global\advance\secno by 1
   \bigpagebreak\noindent
   {\bf \the\parno.\the\secno.~}}


\define\sectionno{
   \global\advance\secno by 1
   \the\parno.\the\secno.
   }


\define\displayno{
   \global\advance\subsecno by 1
   \the\parno.\the\secno.\the\subsecno
   }


\define\formulano{\tag\displayno}

\define\actualsecno{(\the\parno.\the\secno)}

\define\actualsubsecno{(\the\parno.\the\secno.\the\subsecno)}

%
\define\label#1{\xdef#1{\actualsecno}}

%
\define\sublabel#1{\xdef#1{\actualsubsecno}}

%
\define\plabel#1{\xdef#1{\actualsecno}{\smc\string#1}}

%
\define\psublabel#1{\xdef#1{\actualsubsecno}{\smc\string#1}}


\topmatter
\title  On the signature character of representations of $p$-adic
general linear groups \endtitle

\rightheadtext{Signature of $p$-adic representations}

\leftheadtext{C. Boyallian, T. Wedhorn}

\author{ Carina Boyallian ,  Torsten Wedhorn} \endauthor

\address CIEM-Famaf, Universidad Nacional de C\'ordoba - (5000) C\'ordoba,
Argentina \endaddress
\email boyallia\@mate.uncor.edu \endemail

\address  Mathematisches Institut der Universität zu K\"oln, Weyertal
86-90, D-50931 K\"oln, Germany  \endaddress
\email  wedhorn\@MI.Uni-Koeln.DE   \endemail

  
\abstract
In this article we calculate the signature character of certain Hermitian
representations of $GL_N(F)$ for a $p$-adic field $F$. We further give a
conjectural description for the signature character of unramified
representations in terms of Kostka numbers.

\endabstract
\endtopmatter


\document

\subheading {0. Introduction}

\smallskip

Let $F$ be a finite extension of the field $\QQ_p$ of $p$-adic rational
numbers and let $G^{\vee}$ be a connected reductive group over $F$. Let $I
\subset G(F)$ be an Iwahori subgroup. We consider irreducible admissible
(complex) representations $\Cal{V}$ of $G(F)$ such that the space of
Iwahori fixed vectors $\Cal{V}^I$ is nonzero. It is well known that those
representations correspond to representations of the Hecke algebra
$\Cal{H}$ associated to the extended Weyl group $\tilde{W}$ of
$G$. Moreover Barbasch and Moy proved in \cite{BM1} that $V$ is Hermitian
(resp.\ unitary) if and only if $\Cal{V}^I$ is Hermitian (resp.\ unitary)
with respect to a certain involution of $\Cal{H}$.

The Hecke algebra $\Cal{H}$ contains the group algebra $\CC[W]$ of the Weyl
group of $G$ as a subalgebra. In particular we can consider for every
irreducible $\Cal{H}$-module $V$ and for every irreducible representation
$\lambda$ of $W$ the space $V_{\lambda} = \Hom_W(\lambda,V)$. If $V$ is
Hermitian, $V_{\lambda}$ inherits a Hermitian form. In this article we
study the signatures $\sigma_{\lambda}(V)$ of this induced Hermitian form
for $G = GL_N$. We call the tuple $(\sigma_{\lambda}(V))_{\lambda}$ the
signature character of $V$.

To do this we reduce the problem to an
analogous one for modules with real central character of a certain graded
Hecke algebra $\HH$ defined by Lusztig \cite{Lu2}. Again we can consider
$\CC[W]$ as a subalgebra of $\HH$. For these $\HH$-modules there is a
classification analogous to the Langlands classification which parametrizes
irreducible representations in terms of subsets $S$ of a chosen basis of
roots, a tempered representation $U$ of the standard Levi subgroup $M$
corresponding to $S$ and a dominant real character $\nu$ of the center of
$M$. These irreducible $\HH$-modules $L(S,U,\nu)$ are quotients of so
called standard modules $X(S,U,\nu)$.

Our strategy to calculate the signature character is the following: We fix
$S$ and $U$ as above. The dominant real characters $\nu$ form a cone. Those
$\nu$ such that $X(S,U,\nu)$ is reducible define affine hyperplanes in this
cone. We call these hyperplanes ``reducibility walls''. For $\nu$ outside
the reducibility walls the signature character is locally
constant in $\nu$. Moreover, ``nearby zero'' Tadi\v c's classification of
unitary representations of $GL_N$ \cite{Ta} implies that $X(S,U,\nu) =
L(S,U,\nu)$ is unitary. As it is possible to determine the
$\CC[W]$-structure of the standard modules, we know the signature character
of $L(S,U,\nu)$ for small $\nu$.

On the other hand a limit argument by Barbasch and Moy \cite{BM3} also
gives an expression of the Hermitian form for $\nu$ ``nearby infinity''
which allows us to express the signature character purely in terms of the
character of the symmetric groups and certain Kostka numbers. For example
our description of the signature at infinity implies that for $S =
\emptyset$ we have $\sigma_{\lambda} = \chi_{\lambda}(w_0)$ at infinity
(where $\chi$ denotes the character of representations of $W$ and $w_0$ is
the longest element in $W$). Hence to calculate the remaining signature
characters we fix a $\nu_0$ lying on a single reducibility wall and a one
paramater family $t \mapsto \nu(t)$ for small $t$ with $\nu(0) = \nu_0$
such that $\nu(t)$ does not lie on the reducibility wall for $t \not=
0$. We then would like to express the sum or the difference of the
signature characters of $L(S,U,\nu(t))$ for positive and negative $t$.

We cannot do this for all reducibility walls. Instead of this we
concentrate on those walls which are needed to calculate the signature
character for unramified representation (i.e.\ for those representations
with $S = \emptyset$). Here we give a precise conjecture for the difference
and the sum of signature characters of both sides. Moreover we prove this
conjecture for unramified representations. As a consequence we get an
explict expression of the signature character of unramified representations
$L(\emptyset,\bold{1},\nu)$ with $L(\emptyset,\bold{1},\nu) =
X(\emptyset,\bold{1},\nu)$ in terms of Kostka numbers. We also give a
description of the signature characters on the reducibility walls for
unramified representations.

\smallskip

We will now give a short overview over the organization of our work: In the
first chapter we describe the equivalence between irreducible representations
$\Cal{V}$ of $G^{\vee}(F)$ with $\Cal{V}^I \not= (0)$ and irreducible
representations with real central character of the associated graded Hecke
algebra $\HH$. The second chapter contains three classifications of
irreducible $\HH$-modules. The first is in terms of conjugacy classes of
pairs $(s,e)$ of a semisimple element $s$ and a nilpotent element $e$ in
the Lie algebra $\frak{gl}_N$. The second one is the Langlands
classification already described above. And the third classification is the
translation of the Bernstein-Zelevinsky classification of representations
of $GL_N(F)$ in terms of supercuspidal representations to the setting of
graded Hecke algebras. We also explain how to obtain one classification
from one of the other ones.

In the third chapter we introduce the Zelevinsky involution which will
allow us to calculate also the signature character of irreducible
representations $L(S,U,\nu)$ which are a proper quotient of the standard
module $X(S,U,\nu)$. The content of the fourth chapter is the description
of the $W$-module structure of the standard modules.

The fifth chapter contains the classification of Hermitian and unitary
$\HH$-modules and the formal definition of the signature character of an
irreducible $\HH$-module. In the sixth chapter we express the Hermitian
form on standard modules in terms of a certain intertwining operator. Here
we follow closely Barbasch and Moy \cite{BM3}. Chapter Seven contains the
description of the reducibility walls and also the theorem that for
unramified representations the isolated unitary representions are precisely
those which lie on the intersection of $[N/2]$ reducibility walls.
In the eighth chapter we give an explicit algorithm to compute the
signature ``nearby infinity''. We do this by making more precise a
description of Barbasch and Moy given in \cite{BM3}.

The nineth chapter now deals with the topic of crossing the reducibility
walls. We define the ``height'' of such a wall and study reducibility walls
of height 1 and 2 in more detail as those are the walls which occur in the
unramified case. Then we prove our wall crossing theorems in the unramified
case. We further formulate a conjecture for the general case of crossing
reducibility walls of height one.

In the tenth chapter we give a conjecture for crossing certain walls of
height bigger than one and use these conjectures to give an explicit
description of the signature character for unramified representations in
terms of Kostka numbers. Finally in the last chapter we calculate the
signature character for all Hermitian representation of $GL_N$ for $N =
2,3$ and $4$.

\medskip

{\bf Notations}: We use the following notations: All algebraic
varieties and all representations are assumed to be over the complex
numbers $\CC$.

If $R$ is any ring, we denote by $M_n(R)$ the ring of $(n \times
n)$-matrices. If $A \in M_n(R)$ and $B \in M_m(R)$ are two matrices we
denote by $A \oplus B \in M_{n+m}(R)$ the bloc matrix $\left(\smallmatrix A
& 0 \\ 0 & B \endsmallmatrix\right)$, and $\text{diag}(\alpha_1,\dots,\alpha_N)
\in M_N(R)$ denotes the diagonal matrix with entries
$\alpha_1,\dots,\alpha_N \in R$. Finally let $I_N \in GL_N(R)$ be the
identity matrix.

If $W$ is a finite group, we denote by $\hat W$ the set of isomorphism
classes of irreducible representations of $W$.


\newparagraph{Representations of p-adic groups and graded Hecke algebras}

\secstart Let $G^{\vee}$ be a split (connected) reductive
group over a $p$-adic field $F$ and let $G$ be the Langlands dual group over
$\CC$. Let $\Cal R = (X,Y,R,R^{\vee},\Pi)$ be the based root datum of $G$
and $s_{\alpha} \in GL(X)$ the reflections associated to the roots $\alpha \in
R$. Denote by $W$ the Weyl group of $\Cal R$, i.e.~ the group generated by the
$s_{\alpha}$ for $\alpha \in R$. The root base $\Pi$ defines a partial order
on $R^{\vee}$ by $\agcheck_1 \leq \agcheck_2$ if $\agcheck_2 - \agcheck_1$ is
a linear combination with nonnegative integer coefficients of elements of
$\{\agcheck\,|\,\alpha \in \Pi\}$. We denote by $\Pi_m$ the set of all
$\beta \in R$ such that $\bgcheck$ is a minimal element for this partial
order.

Elements in the semidirect product $\Wtilde =
W \ltimes X$ are written in the form $wa^x$ for $w \in W$ and $x \in
X$. We set
$$\Stilde = \{s_{\alpha}\,|\, \alpha \in \Pi\} \cup
\{s_{\alpha}a^{\alpha}\,|\, \alpha \in \Pi_m\} \subset \Wtilde.$$
The root base $\Pi$ defines a system of positive roots $R^+ \subset R$
and a length function
$$\eqalign{l\colon \Wtilde &\arr \NN_0, \cr
l(wa^x) &= \sum_{\alpha \in R^+ \atop w(\alpha) \in R^-} \vert \langle
x,\agcheck \rangle + 1\vert + \sum_{\alpha \in R^+ \atop w(\alpha) \in R^+}
\vert \langle x,\agcheck \rangle \vert\cr}$$
which extends the usual length function on the Coxeter subgroup of
$\Wtilde$ which is generated by $\Stilde$.

We further set
$$R^{\flat} = \{\alpha \in R\,|\, \alpha^{\vee} \in 2Y\}.$$
For example if $(X,Y,R,R^{\vee},\Pi)$ is simple and simply connected
(i.e. it is indecomposable and $R^{\vee}$ generates $Y$), $R^{\flat}$ is
nonempty iff the Dynkin type is $B_n$.

\secstart\label{\padictoaffine} We are interested in the following category
of representations of $G^{\vee}(F)$: Let $I \subset G^{\vee}(F)$ be an Iwahori
subgroup. We call a smooth representation $V$ on a complex vector space
$I$-spherical if it is of finite length and if every subquotient of $V$ is
generated by its fixed vectors with respect to $I$.

Denote by $\Hscr(G^{\vee}//I)$ the Hecke algebra of $G^{\vee}(F)$ with
respect to $I$. The underlying vector space consists of the $\CC$-valued
functions of the (discrete) quotient $I\backslash G^{\vee}(F)/I$ with finite
support, and the algebra structure is given by convolution. For every
smooth representation $V$ of $G^{\vee}(F)$ the space of $I$-fixed vectors
$V^I$ is naturally an $\Hscr(G^{\vee}//I)$-module and the functor $V \mapsto
V^I$ induces an equivalence between the category of $I$-spherical
representations of $G^{\vee}(F)$ and the category of left
$\Hscr(G^{\vee}//I)$-modules which are of finite length (or equivalently
finite dimensional as $\CC$-vector spaces).

\secstart\label{\defineaffine} The Hecke algebra $\Hscr(G^{\vee}//I)$ can be
described directly in terms of generators and relations using the based
root datum $\Cal R$. More precisely, it can be considered as a
specialization of the affine Hecke algebra $\Hscr = \Hscr_{\Cal R}$
associated to $\Cal R$ which is defined as follows:

Let $B$ be the braid group of $\Cal R$, i.e. the group with generators $T_w$
for $w \in \Wtilde$ and relations $T_wT_{w'} = T_{ww'}$ whenever $l(w) +
l(w') = l(ww')$. Denote by $z$ an indeterminate. Then $\Hscr$ is the
$\CC[z,z^{-1}]$-algebra which is the quotient of the group algebra (over
$\CC[z,z^{-1}]$) of the braid group $B$ by the two sided ideal generated by
the elements
$$(T_s + 1)(T_s - z^2)$$
for $s \in \Stilde$. For $w \in \Wtilde$ we denote the image of $T_w$ in
$\Hscr$ again by $T_w$.

Let $q$ be the number of elements in the residue field (of the ring of
integers) of $F$ and denote by $\xi_q\colon \CC[z,z^{-1}] \arr \CC$ the
$\CC$-algebra homomorpism which sends $z$ to $q$. Then a classical result of
Iwahori and Matsumoto [IM] shows that $\Hscr(G^{\vee}//I)$ is isomorphic to
$\Hscr \otimes_{\CC[z,z^{-1}],\xi_q} \CC$.

\secstart\label{\affineprop} We collect some properties of $\Hscr$: Let
$$X_{\text{dom}} = \{x \in X\,|\,\hbox{$\langle x, \agcheck \rangle \geq 0$
for all $\alpha \in \Pi$}\}$$
be the set of dominant weights. We define for $x \in X$ an element $\Tbar_x
\in B$ as follows: Write $x = x_1 - x_2$ with $x_1,x_2 \in X_{\text{dom}}$
and set
$$\Tbar_x = T_{a^{x_1}}T_{a^{x_2}}^{-1}.$$
This is easily seen to be well-defined. Its image in $\Hscr$ is again
denoted by $\Tbar_x$. Further set
$$\theta_x = z^{-l(\Tbar_x)}\Tbar_x \in \Hscr$$
where $l\colon B \arr \ZZ$ is the unique extension of the length function
$l$ to $B$.

Let $\Cal O$ be the group algebra of $X$ over the ring $\CC[z,z^{-1}]$. Then
the map $x \asr \theta_x$ defines an embedding of
$\CC[z,z^{-1}]$-algebras $\Cal O \arr \Hscr$ (\cite{Lu2} 3.4). In the sequel
we consider $\Cal O$ as a subring of $\Hscr$. We further have (loc.\ cit.\
3.7 and 3.11):

\proclaim{Proposition} The Hecke algebra $\Hscr$ is a free left $\Cal
O$-module and a free right $\Cal O$-module with basis $\{T_w\mid w\in W\}$
in each case. The center $\Cal Z$ of the Hecke algebra $\Hscr$ consists the
$W$-invariants of $\Cal O$.
\endproclaim

\secstart Every finite-dimensional left $\Hscr$-module $V$ admits a
(unique) primary decomposition with respect to the center $\Cal Z$ of $\Hscr$
$$V = \bigoplus_{\chi} V_{\chi}$$
where $\chi$ runs through the set of characters of the center $\Cal Z$. As
$\Cal Z = \Cal O^W = \CC[z,z^{-1}][X]^W$, the characters of $\Cal Z$ are given
by pairs $(Wt,z_0)$ where $Wt$ is a $W$-orbit of $T = Y \otimes_{\ZZ}
\CC^\times$ and where $z_0 \in \CC^\times$. This induces a decomposition of the
category of finite-dimensional left $\Hscr$-modules into the direct sum of
the categories ${}_{\Hscr}\Mod_{\chi}$ of finite-dimensional left
$\Hscr$-modules $V$ such that $\zeta - \chi(\zeta)$ is nilpotent on $V$
for all $\zeta \in \Cal Z$. If $\chi$ corresponds to the pair $(Wt,z_0)$, we
also write ${}_{\Hscr}\Mod_{Wt,z_0}$.

\secstart\label{\definegraded} The categories ${}_{\Hscr}\Mod_{Wt,z_0}$
are equivalent to categories of finite-dimensional left $\HH$-modules where
$\HH$ is a graded version of $\Hscr$. Instead of explaining how to obtain
$\HH$ from $\Hscr$ by grading with respect to a certain ideal, we give the
abstract definition of the graded Hecke algebra $\HH = \HH_{\Cal R}$
associated to a reduced root datum $\Cal R = (X,Y,R,R^{\vee},\Pi)$.

We set $\OO = \CC[r] \otimes_{\CC} \Sym(X \otimes_{\ZZ} \CC)$. It carries
an action by $W$ induced by the trivial action on $\CC[r]$ and the
canonical action on $X$. As a $\CC$-vector space we have
$$\HH = \OO \otimes_{\CC} \CC[W]$$
with a structure of associative
$\CC$-algebra with unit $1 \otimes e$, defined by the rules:
\roster
\item"{(i)}" $\OO \arr \HH$, $\theta \asr \theta \otimes e$ is an
algebra homomorphism.
\item"{(ii)}" $\CC[W] \arr \HH$, $w \asr 1 \otimes w$ is an algebra
homomorphism.
\item"{(iii)}" $(\theta \otimes e)(1 \otimes w) = \theta \otimes w$ for
$\theta \in \OO$ and $w \in W$.
\item"{(iv)}" For all $\alpha \in \Pi$ and $\theta \in X$ we have
$$(1 \otimes s_{\alpha})(\theta \otimes e) - ({}^{s_{\alpha}}\theta \otimes
e)(1 \otimes s_{\alpha}) = 2r\alpha\langle \agcheck, \theta \rangle.$$
\item"{(v)}" $r$ is in the center of $\HH$.
\endroster

Usually we will omit the $\otimes$ when denoting elements in
$\HH$ and write $t_w$ instead of $1 \otimes w$. For $\alpha \in R$ we
further set $t_{\alpha} := t_{s_{\alpha}}$. For example the relation (iv)
becomes
$$t_{\alpha}\cdot\theta - {}^{s_{\alpha}}\theta\cdot t_{\alpha} = 2r\langle
\agcheck, \theta \rangle.$$
Note that this relation is equivalent to
$$\theta\cdot t_{\alpha} - t_{\alpha}\cdot{}^{s_{\alpha}}\theta = 2r\langle
\agcheck, \theta \rangle.$$

\secstart\label{\centralchar} The center of $\HH$ consists of the
$W$-invariants of $\OO$. As above we get a decomposition of the category of
finite-dimensional left $\HH$-modules into categories
${}_{\HH}\Mod_{Wu,r_0}$ where $(Wu,r_0)$ runs through the set of
pairs consisting of a $W$-orbit of an elements $u \in Y \otimes_{\ZZ} \CC$
and a complex number $r_0$.

In the sequel we will consider a character of
the center of $\HH$ also as a pair $(\{s\},r_0)$ where $\{s\}$ is a
$G$-conjugacy class of a semisimple element of $\frak g$ and where $r_0$ is a
complex number.

\secstart\label{\affinegradedsetup} The relation between $\Hscr$-modules
and $\HH$-modules is the following (see \cite{Lu2} 8-10 and
\cite{Lu3} 4): As we are mainly interested in the case $G = GL_N$ we make
the following additional assumptions to simplify notations:

\roster
\item"{$\bullet$}" The derived group of $G$ is simply connected (or
equivalently the center of $G^{\vee}$ is connected).
\item"{$\bullet$}" $R^{\flat} = \emptyset$.
\endroster

Fix a central character of $\Hscr$ corresponding to a pair $(Wt,z_0)$. To
simplify notations further we assume that $z_0$ is a positive real
number, different from 1. Further we choose an auxiliary $t$ in the Weyl
orbit $Wt$. The following constructions will be independent of this choice
(up to isomorphism which is given by $w \in W$ if $t$ is replaced by
$wt$). We can decompose $T = \CC^\times \otimes_{\ZZ} Y$ into an elliptic and a
hyperbolic part, namely $T = T_{\text{ell}} \times T_{\text{h}}$ where
$$T_{\text{ell}} = \{z \in \CC^\times\,|\,|z| = 1\} \otimes Y, \qquad
T_{\text{h}} = \RR_{>0} \otimes Y.$$
This is a decomposition of real Lie groups. We can therefore write uniquely
$t = t_et_h$ where $t_e \in T_{\text{ell}}$ and $t_h \in T_{\text{h}}$.

Let $\Cal R(t) = (X,Y,R(t),R^{\vee}(t),\Pi(t))$ be the root datum with
$$\eqalign{R(t) &= \{\alpha \in R\,|\,\alpha(t_e) = 1\}, \cr
R^{\vee}(t) &= \{\agcheck \in R^{\vee}\,|\,\alpha \in R(t)\}.\cr}$$
Then $R^+(t) = R(t) \cap R^+$ is a system of positive roots and the
associated set of simple roots is $\Pi(t)$. We write $W(t)$ for the
Weyl group of the root datum $\Cal R(t)$.

Denote by
$$\log_{z^2_0}\: T_h = \RR_{>0} \otimes{\ZZ} Y \arr \RR \otimes{\ZZ} Y$$
the isomorphism induced by
$$\RR_{>0} \arriso \RR,\qquad x \asr \log(x)/2\log(z_0).$$

Note that the isomorphism class of $\Cal R(t)$ does not depend on the choice
of $t$ in its Weyl orbit

The relation between $\Hscr$-modules and $\HH$-modules is now given by the
next proposition.

\proclaim{\sectionno Proposition}\label{\affinegraded} The categories
 ${}_{\Hscr_{\Cal R}}\Mod_{Wt,z_0}$ and ${}_{\HH_{\Cal
 R(t)}}\Mod_{W(t)\log_{z^2_0}(t_h),1/2}$ are equivalent. More precisely,
 denote by $\hat\Cal H$ (resp.\ $\hat\HH$) the completion of $\Cal H_{\Cal R}$
 (resp.\ $\HH_{\Cal R(t)}$) with respect to the maximal ideal of the center
 corresponding to $(Wt,z_0)$ (resp.\ $(W(t)\log_{z^2_0}(t_h),1/2)$) (hence
 ${}_{\Hscr_{\Cal R}}\Mod_{Wt,z_0}$ resp.\ ${}_{\HH_{\Cal
 R(t)}}\Mod_{W(t)\log_{z^2_0}(t_h),1/2}$ is the category of modules of finite
 length over $\hat\Cal H$ resp.\ $\hat\HH$). Then
 $\hat\Cal H $ and $\hat\HH$ are Morita equivalent.
\endproclaim

\demo{Proof} This follows from the main result of \cite{Lu2} 8-10 and
\cite{Lu3} 4 with simpler notation due to our assumptions. We also replace
the indeterminate $r$ in loc. cit.\ by $r/2\log(z_0)$. First of all we can
replace $\Cal R$ by its derived root datum $\text{der}(\Cal R)
= (X_{\text{der}}, Y_{\text{der}}, R_{\text{der}}, R^{\vee}_{\text{der}},
\Pi_{\text{der}})$, i.e. $Y_{\text{der}} = Y \cap \langle \alpha \mid
\alpha \in R \rangle_{\RR}$. Then $Y_{\text{der}}$ is a direct summand of
$Y$ and if we set $T_{\text{der}} = Y_{\text{der}} \otimes \CC^\times$, we
get a decomposition $T = T' \times T_{\text{der}}$ such that $\alpha(t') =
1$ for all $\alpha \in R$ and hence a decomposition $t =
t't_{\text{der}}$. Now $R(t)$, $R^{\vee}(t)$ and $W(t)$ depend only on the
derived root datum and on $t_{\text{der}}$.

Now our assumption in \affinegradedsetup\ implies that we are in the
situation of \cite{Lu3} 4.4. Then loc.\ cit.\ 4.5 shows that the group
$\Gamma_t$ defined in \cite{Lu2} 8.1 is trivial which implies the result as
explained in \cite{Lu3} 4.9.\qed\enddemo

\secstart If $\Cal R$ is the root datum of $G = GL_N$ or more generally of a
reductive group $G$ whose simple components are of Dynkin type $A_n$, then for
any choice of $t$ the root datum $\Cal R(t)$ is the root datum of a Levi
subgroup of $G$, namely of the centralizer of $t_e$ (after considering $T$ as
a maximal torus of $G$ which is well defined up to conjugacy). For other
Dynkin types this is not true in general.


\newparagraph{Classification of $\HH$-modules}

\secstart\label{\simplesetup} From now on we make the simplifying
assumption that $G$ is isomorphic to a Levi subgroup
of $GL_N$ (hence isomorphic to a product of
groups of the form $GL_{N_i}$). Let $\Cal R = (X,Y,R,R^{\vee},\Pi)$ the
based root datum of $G$, $W$ its Weyl group, and denote by $\HH$ the
associated Hecke algebra. We fix $r_0\in \CC \setminus \{0\}$, and denote by
$\Irr_{r_0}(\HH)$ the set of equivalence classes of irreducible
$\HH$-modules, where $r \in \HH$ acts by $r_0$.

\secstart\label{\definecenter} We set ${\frak t} = Y \otimes \CC$ and
denote by $\frak t^* = X \otimes \CC$ the dual. Define further
$$\align
\frak z\phantom{{}^*} &= \{\lambda \in \frak t\phantom{{}^*} \mid
\text{$\langle \lambda,\alpha \rangle = 0$ for all $\alpha \in R$}\},\\
\frak z^* &= \{x \in \frak t^* \mid \text{$\langle \agcheck,x \rangle = 0$
for all $\agcheck \in R^{\vee}$}\}.
\endalign$$
Then $\frak z$ is canonically isomorphic to the center of $\Lie(G)$ and the
duality of $\frak t$ and $\frak t^*$ induces a perfect duality of $\frak z$
and $\frak z^*$. The subspace $\frak z^*$ of $X \otimes \CC$ has as a
complement the space generated by $\alpha \in R$.

Every element $\nu \in {\frak z}$ defines a one-dimensional $\HH_G$-module
$\CC_{\nu} = \CC_{\nu,r_0}$ where $r$ acts by multiplication with $r_0$,
$\alpha \in R$ acts by multiplication with $2r$, $\xi \in \frak z^*$ acts
by $\nu(\xi)$ and $s_{\alpha}$ acts trivially.

\secstart\label{\Levisubalgebra} Let $S \subset \Pi$ be a subset and denote by
$(X,Y,R_S,R_S^{\vee},S)$ the corresponding subroot system with Weyl group
$W_S$ and associated graded Hecke algebra $\HH_S$. For this root system we
have the subspaces $\frak z_S$ and $\frak z^*_S$ as in \definecenter. Note
that $\HH_S$ is isomorphic to the graded Hecke algebra associated to a Levi
subgroup $M$ of $G$ and that $M$ again satisfies the condition in \simplesetup.

We have a canonical embedding $\HH_S \subset
\HH$ which makes $\HH$ into a free $\HH_S$-module of rank $\#(W/W_S)$.
If $U$ is a $\HH_S$-module we also write $\Ind^{\HH}_{\HH_S}(U)$ instead of
$\HH \otimes_{\HH_S} U$.

\secstart\label{\KLclass} By \cite{Lu1I} there exists a bijection between
$\Irr_{r_0}(\HH)$ and the set of $G$-conjugacy classes of pairs of the
form $(s,e)$ where $s\in\frak g$ is semisimple, $e\in\frak g$ is nilpotent
such that
$$
[s,e]=2r_0e.
$$
Denote by $L_G(s,e,r_0) = L(s,e)$ the irreducible $\HH$-module
corresponding to the conjugacy class $\{(s,e)\}$ of $(s,e)$. By
\cite{Lu5} 1.15, $L(s,e)$ is the unique irreducible quotient of a
standard module $X_G(s,e,r_0) = X(s,e)$ associated to $\{(s,e)\}$. See
\cite{Lu1} 8 and \cite{Lu4} 10 for the definition of the $X(s,e)$
(and \cite{Lu4} 10.11 for the fact that both constructed modules are
isomorphic).

\secstart By definition of the standard module $X_G(s,e,r_0)$
\cite{Lu1} 8, its central character exists and is given by the pair
$(\{s\},r_0)$ \centralchar. In particular this is also true for the
irreducible quotient $L_{r_0}(s,e)$.

\secstart\label{\Ldata} Next we will describe the Langlands
classification. Let $s \in {\frak g} = \frak{gl}_N$ be semisimple
and $e \in \frak{gl}_N$ nilpotent such that $[s,e]=2r_0e$. First choose a
homomorphism of Lie algebras
$$
\psi\: sl_2\longrightarrow \frak{gl}_N
$$
such that $e = \psi \left(\smallmatrix 0 & 1\\ 0 & 0
\endsmallmatrix\right)$ and such that $[s,f] = -2r_0f$ where $f =
\left(\smallmatrix 0 & 0\\ 1 & 0 \endsmallmatrix\right)$. This implies
$[s,h] = [s,[e,f]] = 0$ by the Jacobi identity if $h = \psi
\left(\smallmatrix 1 & \phantom{-}0\\ 0 & -1 \endsmallmatrix\right)$. Such a
homomorphism exists and is uniquely determined up to conjugation with $z
\in G$ such that $\text{Ad}(z)e = e$ and $\text{Ad}(z)s = s$ by a variant
of the Jacobson-Morozov theorem (see \cite{KL} 2.4(g) and 2.4(h)).

\secstart\label{\definetempered} From now on assume that $r_0$ is a
positive real number and that the conjugacy class of $s$ is in $(Y
\otimes_{\ZZ} \RR)/W$ and define $t := s - r_0h$. Then $t$ is a
semisimple element of $\frak g$ whose conjugacy class depends only on the
conjugacy class of $(s,e)$. Note that $t$ commutes with $s$ and with the
image of $\psi$. We call the conjugacy class of $(s,e)$ (or the associated
simple module $L(s,e)$) {\it tempered} if $t = 0$.

If $\{(s,e)\}$ is tempered, we have $L(s,e) = X(s,e)$ by \cite{Lu5}
1.21.

\secstart\label{\Langlandsclass} We state the following version of
Langlands classification for $\HH$-modules which is a slight reformulation
of \cite{Ev} (recall that we assume that $r_0 \in \RR$ and $\{s\} \in (Y
\otimes_{\ZZ} \RR)/W$):

\proclaim{Theorem} For every irreducible $\HH$-module $V$ there exists
a triple $(S,U,\nu)$ where $S$ is a subset of $\Pi$, $U$ is an irreducible
tempered representation of $\HH_S$ and where
$$
\nu \in \frak z^+_S = \{\lambda \in \frak z_S \cap Y \otimes \RR \mid
\text{$\langle \lambda, \alpha \rangle > 0$ for all $\alpha \in \Pi
\setminus S$}\}
$$
such that $V$ is the unique irreducible quotient of $\Ind^{\HH}_{\HH_S}(U
\otimes \CC_{\nu})$. Further $S$ and $\nu$ and the isomorphism class of $U$
are uniquely determined by the isomorphism class of $V$. We set $L(S,U,\nu)
= V$ and $X(S,U,\nu) = \Ind^{\HH}_{\HH_S}(U \otimes \CC_{\nu})$.

If $L(S',U',\nu')$ is any other irreducible subquotient of
$\Ind^{\HH}_{\HH_S}(U \otimes \CC_{\nu})$ coming from a triple
$(S',U',\nu')$, we have $\nu' < \nu$ (with respect to the extended Bruhat
order on $\Wtilde$).
\endproclaim

The triple $(S,U,\nu)$ is called {\it Langlands data} associated to $V$.

\secstart\label{\LtoLanglands} We keep the notations of \Ldata. In
particular we have the $G$-conjugacy class $\{(s,e)\}$ and the associated
irreducible $\HH$-module $L(s,e)$. We are now going to explain how
to obtain the Langlands data $(S,U,\nu)$ corresponding to $L(s,e)$. We
follow \cite{Lu5} 3.9ff choosing for $\tau$ in loc.~cit.\ the
homomorphism $\CC \arr \RR$ which associates to each complex number its
real part. We use a somewhat more explicit but less canonical
description. For this we fix a Borel subgroup $B$ of $G$ and a
maximal torus $T$ of $G$ contained in $G$ which gives an
identification of the abstract based root datum $\Cal R$ with the based root
datum of $(G,B,T)$. In particular we have an identification $\Lie(T) =
\frak t = Y \otimes \CC$. The conjugacy class of $s$ is then an element in
$(Y \otimes \RR)/W$ (because we assumed that $s$ is real). After
$G$-conjugation we can assume that $s$ is a diagonal matrix and that $(s,e)
= (s_1,e_1) \oplus \dots \oplus (s_l,e_l)$ where $e_i$ is of the form
$$
e_i = \pmatrix 0 \\ 1 & 0 \\ & 1 & 0 \\ & & \ddots & \ddots \\ & & & 1 & 0
\endpmatrix.
$$
Let $\sigma_i$ be the first entry in the diagonal matrix $s_i$. We can
assume that $\sigma_1 + \frac12 m_1 \geq \dots \geq \sigma_m + \frac12
m_l$. Now we can choose $\psi$ as in \Ldata\ such that $h \in Y \otimes \RR
\subset \frak t$, in particular $t = s - r_0h \in Y \otimes \RR$. We have
$\langle t, \alpha \rangle \geq 0$ for all $\alpha \in \Pi$.

We set
$$
S = \{\alpha \in \Pi \mid \langle t, \alpha \rangle = 0\}.
$$
Then $\HH_S$ is isomorphic to the graded Hecke algebra associated to Levi
subalgebra $\frak m = \text{Cent}_{\frak g}(t)$. Note that $s$ and the
image of $\psi$ are contained $\frak m$. Let $M$ be the
corresponding Levi subgroup of $G$. As the derived group of $G$ is simply
connected, we have $M = Z_G(t)$. Further we have $[r_0h,e] =
2r_0e$. Hence the conjugacy class of $(r_0h,e)$ defines an irreducible
representation
$$
U := L_M(r_0h,e)
$$
of $\HH_S$. By definition this is a tempered representation of
$\HH_S$. Finally let $\nu$ by the dominant representative of the $W$-orbit
of $t$. By definition we have $\nu \in \frak z^+_S$. We then have $U
\otimes \CC_{\nu} = L_M(s,e) = X_M(s,e)$ and by \cite{Lu5} 3.38 we
have $X_G(s,e) = \Ind_{\HH_M}^{\HH_G}(U \otimes \CC_{\nu})$.

\secstart For the rest of chapter 2 let us assume that $G = GL_N$. Hence
the $p$-adic group $G^{\vee}$ is isomorphic to $GL_N$ and we have the
classification of representations of the $p$-adic group $GL_N(F)$ by
Bernstein-Zelevinsky in terms of supercuspidal representations. By a
theorem of Casselman (see e.g.\ \cite{Ca}
3.8) an irreducible admissible representation of $GL_N(F)$ admits
nontrivial fixed vector under an Iwahori subgroup if and only if its
supercuspidal support consists of unramified quasi-characters. Using
\affinegraded\ we obtain a classification of irreducible $\HH$-modules
where $r \in \HH$ acts by $1/2$ (or equivalently a classification of
irreducible modules of the algebra $\HH^{1/2} = (\HH
\otimes_{\CC[r],r\mapsto 1/2} \CC)$) which we call the BZ-classification. We
further can assume that the central character is real. It is described in
the next sections.

\secstart For $G = GL_1$, $\HH^{1/2}_G$ is nothing but a polynomial algebra in
one indeterminate and we consider any complex number as a one-dimensional
representation of $\HH^{1/2}_{GL_1}$.

Let $\vec{N} = (N_1,\dots,N_m)$ be a tuple of positive integers and set
$GL_{\vec N} = GL_{N_1} \times \dots GL_{N_m}$. If $V_i$
is an $\HH^{1/2}_{GL_{N_i}}$-module ($i = 1,\dots,m$), we set $V_1 \boxtimes
\cdots \boxtimes V_m$ for the module of
$$
\HH^{1/2}_{GL_{\vec N}} = \HH^{1/2}_{GL_{N_1}} \otimes
\cdots \otimes \HH^{1/2}_{GL_{N_m}}
$$
whose underlying vector space is $V_1 \otimes \cdots \otimes V_m$ and which
is endowed with the componentwise action. We further set
$$
V_1 \boxdot \cdots \boxdot V_m = \HH^{1/2}_{GL_N}
\otimes_{\HH^{1/2}_{GL_{\vec N}}} (V_1 \boxtimes \cdots \boxtimes V_m)$$
where $N = N_1 + \cdots + N_m$.

If in particular $(\sigma_1,\dots,\sigma_m)$ is a tuple of complex
numbers, $\sigma_1 \boxdot \dots \boxdot \sigma_m$ is a
$\HH^{1/2}_{GL_m}$-module.

\secstart\label{\definesegment} Given $m \in \NN$ and $\sigma$ a real
number, define the {\it segment}
$$
\Delta(\sigma,m) = [\sigma,\ \sigma + 1 , \ldots, \sigma + m-1].$$
The real number $\sigma + \frac{m-1}{2}$ is called the {\it center} of
$\Delta(\sigma,m)$, and the integer $m \geq 1$ is called the {\it length}
of $\Delta(\sigma,m)$.

Consider $\Delta_1 = \Delta(\sigma_1,m_1)$ and $\Delta_2 =
\Delta(\sigma_2,m_2)$ two segments. We will say that $\Delta_1$ and
$\Delta_2$ are {\it linked} if $\Delta_1\not\subseteq\Delta_2$ and
$\Delta_2\not\subseteq\Delta_1$ and $\Delta_1\cup\Delta_2$ is of the form
$\Delta(\tau,m')$, for some $\tau \in \{\sigma_1,\sigma_2\}$.

Further, we say that $\Delta_1$ {\it precedes} $\Delta_2$, if $\Delta_1$
and $\Delta_2$ are linked and $\tau = \sigma_1$.

\secstart\label{\BZclass} With these definitions we can establish the
following facts (\cite{Ze},\cite{KL}, \cite{Ku}):

\roster
\item Take $\Delta=\Delta(\sigma,m)$ as above. Then,
$\sigma \boxdot (\sigma+1) \boxdot \cdots \boxdot (\sigma + m-1)$ is
reducible for $m > 1$ and has a unique irreducible quotient $L(\Delta)$.

\item Let $(\Delta_1,\ldots,\Delta_l)$ be a tuple of segments as above, and
assume that $\Delta_i$ does not preceed $\Delta_j$ for $i<j$, then
$L(\Delta_1) \boxdot \cdots \boxdot L(\Delta_l)$ admits a unique irreducible
quotient $L(\Delta_1,\ldots,\Delta_l)$.

\item Every irreducible admissible representation of $\HH^{1/2}_{GL_N}$ is
isomorphic to some $L(\Delta_1,\ldots,\Delta_l)$ where
$(\Delta_1,\ldots,\Delta_l)$ is a tuple as in (2). If
$(\Delta'_1,\ldots,\Delta'_k)$ is any other tuple as in (2) such that
$L(\Delta_1,\ldots,\Delta_l) \cong L(\Delta'_1,\ldots,\Delta'_k)$ then $l =
k$ and $\Delta'_i = \Delta_{\pi(i)}$ for some permutation $\pi \in S_l$.

\item $L(\Delta_1) \boxdot \cdots \boxdot L(\Delta_l)$ is irreducible if and
only if no two segments $\Delta_i$ and $\Delta_j$ are linked. 

\endroster

We set
$$
X(\Delta_1,\dots,\Delta_l) := L(\Delta_1) \boxdot \cdots \boxdot L(\Delta_l).
$$

\secstart\label{\KLBZ} Let us connect the BZ-classification with the
classification by conjugacy classes of pairs $(s,e)$ such that $[s,e] = e$
(note that to simplify we are still in the case $r_0 = 1/2$ which we can
assume anyway because of \affinegraded).

Let $V = L(\Delta(\sigma_1,m_1),\ldots,\Delta(\sigma_l,m_l))$ be an
irreducible representation of $\HH^{1/2}_{GL_N}$ with $N = m_1 + \cdots
m_l$. Denote by $\gamma_i = \sigma_i + \frac12 (m_i-1)$ the center of
$\Delta(\sigma_i,m_i)$. We assume that $\gamma_1 \geq \gamma_2 \geq \dots
\geq \gamma_l$ which in particular implies that $\Delta(\sigma_i,m_i)$ does
not precede $\Delta(\sigma_j,m_j)$ for $i < j$. We set
$$
s = \bigoplus_{i=1}^l\text{diag}(\sigma_i,\sigma_i+1,\dots,\sigma_i+m_i-1),
\qquad e = \bigoplus_{i=1}^l n_{m_i}
$$
where $n_d$ is the nilpotent $(d \times d)$-matrix
$$
n_d = \pmatrix 0 \\ 1 & 0 \\ & 1 & 0 \\ & & \ddots & \ddots \\ & & & 1 & 0
\endpmatrix.
$$
Then we have $[s,e] = e$ and the irreducible $\HH^{1/2}_{GL_N}$-module
associated to the conjugacy class of $(s,e)$ is isomorphic to $V$.

Now we can use \LtoLanglands\ to compute the corresponding Langlands triple
$(S,U,\nu)$. For this we have to construct the element $t$: We assume that
the based root datum is given by the Borel pair $T \subset B$ of $GL_N$
where $T$ is the diagonal torus and $B$ the Borel subgroup of upper
triangular matrices. The simple roots in $\Pi$ are then given by the linear
forms $\alpha_i\: \text{diag}(x_1,\dots,x_N) \mapsto x_i - x_{i+1}$ for $i
= 1,\dots,N-1$.

As homomorphism $\psi\colon \frak s\frak l_2 \to \frak{gl}_N$ we choose the
unique $\psi$ such that $\psi\left(\smallmatrix 0 & 1 \\ 0 & 0
\endsmallmatrix\right) = e$ and such that
$\psi\left(\smallmatrix 1 & \phantom{-}0 \\ 0 & -1 \endsmallmatrix\right) =
\text{diag}(h_1,\dots,h_l)$ where $h_i =
\text{diag}(-m_i+1,-m_i+3,\dots,m_i-1)$. Hence we have
$$
t = s - \frac12 h = \bigoplus_{i=1}^l
\text{diag}(\undersetbrace \text{$m_i$ times}\to{\gamma_i,\dots,\gamma_i}).
$$

\secstart We give an example: Let $V$ be the irreducible
$\HH^{1/2}_{GL_6}$-module given by the sequence of segments
$$
([2,3],[0,1,2],[1]).
$$
Then we have $V = L(s,e)$ where $s = \text{diag}(2,3,0,1,2,1)$ and
$$
e = \pmatrix
0 \\
1 & 0 \\
  &   & 0 \\
  &   & 1 & 0 \\
  &   &   & 1 & 0 \\
  &   &   &   &   & 0
\endpmatrix.
$$
Further $t = \text{diag}(5/2, 5/2, 1, 1, 1, 1)$ and hence $S = \Pi
\setminus \{\alpha_2\}$. The tempered representation $U$ can be considered
as an $(\HH^{1/2}_{GL_2} \otimes \HH^{1/2}_{GL_4})$-module and it is the tensor
product of the $\HH^{1/2}_{GL_2}$-module $U_1$ and the
$\HH^{1/2}_{GL_4}$-module $U_2$ where $U_1$ is given (with respect to
the BZ-classification) by $[-1/2,1/2]$ and $U_2$ is given by
$([-1,0,1],[0])$. We have $\frak z_S^+ = \{(x_1,x_1,x_2,x_2,x_2,x_2) \in
\RR^6 \mid x_1 > x_2\}$ and $\nu = (5/2,5/2,1,1,1,1)$.

\secstart\label{\KLirred} We remark that we can also check the
irreducibility of $X(s,e)$ directly: The Levi subgroup $Z_G(s)$ acts on the
vector space $\{n \in \frak{gl}_N \mid [s,n] = n\}$ by conjugation and
$X(s,e)$ is irreducible if and only if $e$ lies in the unique open orbit of
that action.

\secstart\label{\subofstandard} Let $\frak{M}$ be a multiset of segments
and let $X(\frak{M})$ be the corresponding standard module. The irreducible
subquotients of $X(\frak{M})$ can be described as follows ([Ze]):

An elementary operation on a multiset $\frak{M}$ is by definition to take two
segments $\Delta_1$ and $\Delta_2$ which are linked from $\frak{M}$ and to
replace them by $\Delta_1 \cup \Delta_2$ and $\Delta_1 \cap \Delta_2$. For two
multisets $\frak{M}_1$ and $\frak{M}_2$ we say that $\frak{M}_1 \preceq
\frak{M}_2$ if $\frak{M}_1$ can be obtained from $\frak{M}_2$ by elementary
operations.

With this definition we have that $L(\frak{M}')$ occurs as irreducible
subquotient of $X(\frak{M})$ if and only if $\frak{M}' \preceq \frak{M}$.

Further, $X(\frak{M})$ has a unique irreducible quotient and a unique
irreducible submodule and the isomorphism classes of both of them occur
with multiplicity one in $X(\frak{M})$.

Zelevinsky shows in loc.\ cit.\ that if all real numbers which
occur in all segments of $\frak{M}$ are pairwise distinct, $L(\frak{M}')$
occurs in $X(\frak{M})$ with multiplicity one for all $\frak{M}' \preceq
\frak{M}$.

\secstart\label{\defineunramified} Let $V = L(s,e) =
L(\Delta_1,\dots,\Delta_l)$ be an irreducible $\HH^{1/2}_{GL_N}$-module
with associated Langlands data $(S,U,\nu)$. Then we call $V$ {\it
unramified} if the following equivalent coniditions are satisfied:
\roster
\item We have $e = 0$ and $s$ is a regular semisimple element.
\item The length of all segments $\Delta_i$ is equal to 1 and their centers are
pairwise different.
\item $S = \emptyset$.
\endroster

\secstart\label{\chartemp} Let $V = L(\Delta_1,\dots,\Delta_l)$ be an
irreducible $\HH^{1/2}_{GL_N}$-module with associated Langlands data
$(S,U,\nu)$. Then the following assertions are equivalent:
\roster
\item $V$ is tempered.
\item $S = \Pi$ and $\nu = 0$.
\item All centers of the segments $\Delta_i$ are equal to zero.
\endroster


\newparagraph{The Zelevinsky involution}

\secstart\label{\defineZel} Let $G$ be a reductive group over $\CC$ with
based root datum $(X,Y,R,R^{\vee},\Pi)$ and let $\HH_G$ be the associated
graded Hecke algebra. We define an involution on $\HH_G$ which we show to be
induced by the Zelevinsky involution on the affine
Hecke algebra $\Cal H_G$ for $G = GL_N$. By abuse of notation we will call
it also Zelevinsky involution and denote it by $\zeta$. It is defined as
$$\alignat2
\zeta(r) &= r, \\
\zeta(t_w) &= (-1)^{l(w)}t_{w_0ww_0} &&\qquad\text{for $w \in W$,}\\
\zeta(\theta) &= {}^{w_0}\theta &&\qquad\text{for $\theta \in X$}
\endalignat
$$
where $w_0$ is the element of maximal length in the Weyl group $W$ of $G$.

It is easy to check that $\zeta$ preserves the relations in \definegraded\
defining the graded Hecke algebra.

For every $\HH_G$-module $V$ (where the $\HH_G$-module structure is given
by a $\CC$-algebra homomorphism $\rho\: \HH_G \to \End(V)$) we write
$\zeta(V)$ for the $\HH_G$-module given by $\rho \circ \zeta$. This
defines an involutive endofunctor of the category of $\HH_G$-modules
${}_{\HH_G}\Mod$.

\secstart Note that from the definition we get the following observation

\proclaim{Remark} Let $V$ be an $\HH_G$-module which admits a central
character $\chi$. Then $\zeta(V)$ also admits a central character $\chi'$,
and we have $\chi = \chi'$.
\endproclaim

\secstart\label{\defineGroth} We want to describe the effect of $\zeta$ on
irreducible $\HH_{GL_N}$-modules given by the BZ-classification. For this
we show that $\zeta$ induces Zelevinsky's involution. To prove this we make the
following definition: Let $\Cal R_N$ be the Grothendieck group of
the category of finite-dimensional $\HH^{1/2}_{GL_N}$-modules with real
central character and set $\Cal R =
\bigoplus_{N \geq 0}\Cal R_N$ where $GL_0$ is by definition the trivial
group (and hence $\Cal R_0 = \ZZ$). The map
$$
\Cal R_{N_1} \times \Cal R_{N_2} \to \Cal R_{N_1+N_2}, \qquad ([V_1],[V_2])
\mapsto [V_1 \boxdot V_2]
$$
makes $\Cal R$ into a graded ring. For every segment $\Delta$ as in
\definesegment\ we have the corresponding irreducible representation
$[L(\Delta)] \in \Cal R$ and the same arguments given for the analogous
statement for representations of $GL_N(F)$ in \cite{Ze} 7 show that this
makes $\Cal R$ into the polynomial algebra over $\ZZ$ in indeterminates
$\Delta$ where $\Delta$ runs through all segments.

\secstart\label{\compareZel} The Zelevinsky involution defines an involutive
automorphism of the graded ring $\Cal R$. We define another involution
which is by definition the unique involutive automorphism $\zeta'$ of
$\Cal R$ such that
$$
\zeta'(L([\sigma,\sigma+1,\dots,\sigma+m-1])) =
L([\sigma+m-1],[\sigma+m-2],\dots,[\sigma])
$$
which is the analog of the automorphism constructed by Zelevinsky in
\cite{Ze} 9.12.

\proclaim{Proposition} The involutions $\zeta$ and $\zeta'$ of $\Cal R$
coincide.
\endproclaim

\demo{Proof} This follows from results of Moeglin and Waldspurger
\cite{MW} I for the affine Hecke algebra using \affinegraded. We remark
that the elements $X_i$ (resp.\ $S_j$) of loc.~cit.\ are those which are called
$\theta_{e_i}$ (resp.\ $T_{s_j}$) in \affineprop\ where $e_i \in \ZZ^N = X$
is the $i$-th standard base vector and $s_j$ is the reflection
corresponding to the base root $(x_1,\dots,x_N) \mapsto x_j -
x_{j+1}$. Further not that for the image $t_{s_j}$ of $S_j$ in the graded Hecke
algebra we have $t_{s_j}^{-1} = t_{s_j}$. Finally note that $q$ in
loc.~cit.\ is equal to $z^2$ in \defineaffine\ and hence that under the
transition to the graded Hecke algebra as described in \affinegraded\ the
factor $q$ becomes 1.
\qed\enddemo

\secstart\label{\calculateZel} Now let $V =
L(\Delta(\sigma_1,m_1),\dots,\Delta(\sigma_l,m_l))$ be an irreducible
$\HH^{1/2}_{GL_N}$-module. We want to explain the effect of $\zeta$ on
$V$. For this we follow \cite{MW} II: Set $\Delta_i = \Delta(\sigma_i,m_i)$
and let $\frak{M}$ be the multiset (i.e.\ the set with multiplicities) of the
segments $\Delta_i$. By loc.~cit.\ $\zeta(V)$ is the irreducible
representation associated to the mutiset $\frak{M}^{\#}$ of segments where
$\frak{M}^{\#}$ is defined as follows:

For $t \in \RR/\ZZ$ write $\frak{M}_t = \{\Delta_i \in \frak{M} \mid
\sigma_i \equiv t \bmod \ZZ\}$. Then $\frak{M}$ is the disjoint union of
the $\frak{M}_t$ where $t$ runs through $\RR/\ZZ$ and we set
$$
\frak{M}^{\#} = \bigcup_{t \in \RR/\ZZ}\frak{M}_t^{\#}.
$$

Hence we will from now on assume that $\sigma_1 \equiv \dots \equiv
\sigma_l \bmod \ZZ$. We introduce a total order on the set of segments by
saying $\Delta(\sigma_1,m_1) \geq \Delta(\sigma_2,m_2)$ if
$$
\sigma_1 > \sigma_2,\qquad\text{or}\qquad\sigma_1 =
\sigma_2\quad\text{and}\quad\sigma_1+m_1 \geq \sigma_2+m_2.
$$
Further, if $\Delta = \Delta(\sigma,m)$, we set
$$
\Delta^{-} = \Delta(\sigma,m-1).
$$

Let $\delta$ be the biggest real number appearing in one of the segments of
$\frak{M}$ and let $\Delta_{i_0}$ be a segment containing $\delta$ which is
maximal with this property. Necessarily we have $\delta = \sigma_{i_0} +
m_{i_0} -1$. Now define inductively integers $i_1,\dots,i_r$:
\roster
\item"{--}" $\Delta_{i_s}$ is a segment of $\frak{M}$ preceding
$\Delta_{i_s-1}$ such that $\sigma_{i_s} + m_{i_s}-1 = \delta-s$ and such that
$\Delta_{i_s}$ is maximal with this property.
\item"{--}" $i_r$ is the last integer which can be defined this way.
\endroster
We set
$$
\frak{M}^- = (\Delta'_1,\dots,\Delta'_l)
$$
where
$$
\Delta'_i = \cases \Delta^-_i,&\text{if $i \in \{i_0,\dots,i_r\}$} \\
\Delta_i,&\text{otherwise.} \endcases
$$
Note that $\Delta'_i$ can be empty. Now we define
$$\frak{M}^{\#} = \{\Delta(\delta-r,r+1)\} \cup (\frak{M}^-)^{\#}$$
and proceed inductively.

\secstart We give an example of the effect of $\zeta$: For
$$\align
\frak{M}\phantom{{}^{\#}} &= ([3,4],[2,3,4],[1,2],[1/2],[0],[-1/2],[-1,0,1])\\
\intertext{we get}
\frak{M}^{\#} &= ([4],[4],[3],[1,2,3],[0,1,2],[0],[-1/2,1/2],[-1]).
\endalign$$

\secstart There has also been given another combinatorial description of
$\zeta$ in \cite{KZ}.


\newparagraph{$W$-structure of standard modules}

\secstart We assume in this chapter that $G = GL_N$ hence we have
$W = S_N$. As $\HH_G$ contains $\CC[W]$ as a subalgebra, every
$\HH_G$-module $V$ has also the structure of a representation of $W$. We
are interested in the $\CC[W]$-module structure of the standard modules
$X(s,e)$ and the irreducible quotient $L(s,e)$.

\secstart Let us briefly recall some facts of the theory of representations
of $W = S_N$. The set of isomorphism classes of irreducible representations
of $S_N$ is denoted by $\hat S_N$. We have two distinguished (irreducible)
representations of $S_N$, namely the trivial representation $\text{\bf 1}$
and the sign representation $\text{sgn}$. These are the only 1-dimensional
representations of $S_N$.

Denote by $\Cal P(N)$ the set of partitions of $N$. Given a
partition $\text{\bf d} = [d_1 \geq \dots \geq d_N \geq 0] \in \Cal P(N)$,
define the transpose of $\text{\bf d}$ as $\text{\bf d}^t=[d_1^t \geq \dots
\geq d_N^t \geq 0]$ with $d^t_i = \#\{j \mid d_j \geq i\}$. For $\text{\bf
d} \in S_N$ we set $S_{\text{\bf d}} = S_{d_1} \times \dots \times S_{d_N}$
which we embed into $S_N$ in the usual way.

The set $\hat S_N$ is in bijection to $\Cal P(N)$: 
The representation $\pi_{\text{\bf d}}$ corresponding to $\text{\bf d} \in
\Cal P(N)$ is the unique irreducible representation of $S_N$ such that

\roster
\item"{(a)}" The restriction of $\pi_{\text{\bf d}}$ to the subgroup
$S_{\text{\bf d}}$ of $S_N$ contains a copy of the trivial
representation.

\item"{(b)}" The restriction to $S_{\text{\bf d}^t}$ contains a copy of
the sign representation.
\endroster

The tensor product with the sign representation defines an involution on
$\hat S_N$. More precisely, we have $\pi_{\text{\bf d}^t} \cong
\pi_{\text{\bf d}} \otimes \text{sgn}$. The dimension of $\pi_{\text{\bf
d}}$ is the number of standard Young tableaus of shape $\text{\bf d}$. 

On the other hand, there exists a bijective correpondance between the set
of $GL_N$-orbits of nilpotent elements in $\frak{gl}_N$ and the set $\Cal
P(N)$ given by the block sizes of the Jordan normal form of the nilpotent
element. Combining these two facts we obtain a bijection between nilpotent
orbits in $\frak{gl}_N$ and $\hat S_N$. 

Via this bijection, the principal nilpotent orbit corresponds to the trivial
representation, and the zero orbit corresponds to the sign representation.

On the set of nilpotent $G$-orbits of $\frak g$ there is a partial order
where we say that $\Cal O \leq \Cal O'$ iff $\Cal O$ is contained in the
closure of $\Cal O'$. This corresponds to a partial order on $\Cal P(N)$
which is given by $\text{\bf d} \leq \text{\bf d}'$ iff $d_1 + \dots + d_k
\leq d'_1 + \dots + d'_k$ for all $k = 1,\dots,N$. Hence we get also a
partial order on $\hat S_N$ such that $\text{\bf 1}$ is the greatest
element and $\text{sgn}$ is the smallest element.

\secstart\label{\standWmodule} Let $e\in\frak{gl}_N$ be a nilpotent
element, and $\Cal B_e$ be the variety of Borel subgroups of $GL_N$ that
contain $e$. The Springer correspondence tells us that $H^*(\Cal B_e) =
H^*(\Cal B_e,\CC)$ carries an action of the Weyl group $W$ such that
$H^{\dim(\Cal B_e)}(\Cal B_e)$ is isomorphic to the irreducible
$W$-representation corresponding to the $G$-orbit of $e$.

We want to describe the $W$-action of the standard module
$X_{GL_N}(s,e,r_0)$. If we let $r_0$ vary, these standard modules are by
definition the fibres of a vector bundle with $W$-action over the affine
line. As representations of a finite group cannot be deformed, the
$W$-structure of $X_{GL_N}(s,e,r_0)$ is independent of $r_0$. Hence we can
assume $r_0 = 0$ and we have an isomorphism of $W$-modules (\cite{Lu4}
10.13)
$$
X(s,e) \cong H^*(\Cal B_e) \otimes \text{sgn}.
$$

On the other hand, if $\text{\bf d}$ is the partition corresponding to the
$GL_N$-orbit of $e$, there is an isomorphism of $W$-modules (e.g.~\cite{CP})
$$
H^*(\Cal B_e) \cong \Ind^{S_N}_{S_{\text{\bf d}}}(\text{\bf 1}).
$$

Finally the multiplicity of $\pi_{\text{\bf d}'}$ in
$\Ind^{S_N}_{S_{\text{\bf d}}}(\text{\bf 1})$ is given by the Kostka
number $K_{\text{\bf d}',\text{\bf d}}$ (see e.g.~\cite{Ma} I,6 for a
definition).

Altogether we obtain:

\proclaim{Proposition} Fix an $r_0 \in \CC$. Let $s$ in $\frak{gl}_N$
be a semisimple element and $e \in \frak{gl}_N$ be a nilpotent element
such that $[s,e] = 2r_0e$. Let $\text{\bf e}$ be the partition corresponding
to the $GL_N$-orbit of $e$. Then the $W$-structure of the standard module
$X(s,e)$ is given by
$$
[X(s,e) : \pi_{\text{\bf d}}] = K_{\text{\bf d}^t,\text{\bf e}},
\qquad\text{\bf d} \in \Cal P(N).
$$
In particular, $X(s,0)$ is isomorphic to $\CC[W]$ as a $W$-module.
\endproclaim

\secstart\label{\ZelonWmodule} We want to compare the underlying $W$-module
structures of an $\HH_{GL_N}$-module $V$ and its image under the Zelevinsky
involution as defined in \defineZel. The involution $\zeta$ on $\HH_{GL_N}$
restricts to an involution on the subalgabra $\CC[W]$ which we denote
again by $\zeta$ and which induces an involutive endofunctor $\zeta$
of the category $\text{\bf Rep}(W)$ of representations
of $W$. Its effect is described by the following result:

\proclaim{Proposition} Let $V$ be a representation of $W$. Then we have
$$
\zeta(V) \cong V \otimes \text{sgn}.
$$
\endproclaim

\demo{Proof} This follows directly from the definitions: The endofunctor on
$\text{\bf Rep}(W)$ given by $\zeta$ is isomorphic to the one given by the
involution $w \asr (-1)^{l(w)}w$ on $\CC[W]$.
\qed\enddemo

\secstart\label{\Wofstandard} Let $X(s,e)$ be a standard module. By
\standWmodule\ the $\lambda \in \What$ corresponding to the dual partition
of the Jordan type of $e$ is the unique maximal $\lambda \in \What$
occuring in $X(s,e)$ and we have $[X(s,e) : \lambda] = 1$.

Moreover it follows from [BM1] that the sum $X'$ of all
$\HH^{1/2}_G$-submodules of $X(s,e)$ which do not contain $\lambda$ is a
maximal $\HH^{1/2}_G$-submodule of $X(s,e)$ and that we have $X/X' = L(s,e)$.


\newparagraph{Hermitian and unitary $\HH$-modules}

\secstart\label{\hermitianinvol} We return briefly to the general
notations of the first chapter. The $\CC$-vector space $X \otimes_{\ZZ}
\CC$ has a conjugation coming from the complex conjugation $\CC$ and this
induces a complex anti-linear algebra involution on $\Sym(X \otimes_{\ZZ} \CC)$
which we denote by $\theta \asr \bar\theta$.
For $w \in W$ we denote by $t_w$ the corresponding element in $\CC[W]
\subset \HH_G$. Let $w_0\in W$ be the longest element.

Define the $*$-operation on $\HH_G$ as follows:
$$
\alignat2
t^*_w &= t_{w^{-1}},&&\qquad\text{for $w\in W$},\\
\theta^* &= (-1)^{\deg
\theta}t_{w_0}({}^{w_0}\bar{\theta})t_{w_0},&&\qquad\text{for $\theta \in
\Sym(X \otimes_{\ZZ} \CC)$},\\
r^* &= r.
\endalignat
$$
It is easy to check that this defines a complex anti-linear involution on
the algebra $\HH_G$. We call a finite-dimensional $\HH$-module $X$ {\it
Hermitian} if there is a non-degenerate Hermitian form
$\langle\,,\,\rangle$ on $X$ such that
$$
\langle H\cdot x_1,x_2 \rangle = \langle x_1,H^*\cdot x_2 \rangle
$$
for $H\in\HH$, and $x_1,\, x_2\in X$. By \cite{BM2} 5 this notion of being
Hermitian corresponds to the obvious one if $X$ comes from an admissible
representation of $G^{\vee}(F)$ be the procedure described in
\padictoaffine\ and \affinegraded.

\secstart\label{\hermitiancrit} Now let us again assume that $G = GL_N$
and let $V$ be an irreducible $\HH^{1/2}_{GL_N}$-module with real central
character. We want to express the property that $V$ is Hermitian in terms
of the Langlands and the Bernstein-Zelevinsky classification.

First let $(S,U,\nu)$ be the Langlands data associated to $V$. Then it
follows from \cite{BM3} 1.5 that $V$ is Hermitian if and only if there
exists a $w \in W$ satisfying
$$\align
w(\nu) &= -\nu,\tag1 \\
w(S) &= S,\tag2 \\
w(U) &\cong U.\tag3
\endalign$$
Because of $(1)$ and $\nu \in \frak z_S^+$ we have necessarily $w
\in w_0W_S$ where $W_S$ is the subgroup of $W$ generated by $s_{\alpha}$
for $\alpha \in S$. If we write $\Pi = \{\alpha_1,\dots,\alpha_{N_1}\}$ as
in \KLBZ, the identity $(2)$ then implies that $\alpha_i \in S$ if and
only if $\alpha_{N-i} \in S$ for all $i = 1,\dots,N-1$.

Now assume that $V = L(\Delta_1,\dots,\Delta_l)$. By [Ta] its Hermitian dual is given
by $L(\Delta_1^h,\dots,\Delta_l^h)$ where
$$[x, x+1, \dots, x+m-1]^h = [-(x+m-1),\dots,-x].$$
In particular we see:

\proclaim{Proposition} Let $V =
L(\Delta_1,\dots,\Delta_l)$ with $\Delta_i = \Delta(\sigma_i,m_i)$ be given as
above. Then $V$ is Hermitian if and only if we can group together the segments
to pairs $\Delta_{i_1}$ and $\Delta_{i_2}$ ($i_1$ not necessarily different
from $i_2$) such that
\roster
\item The center of $\Delta_{i_1}$ is the negative of the center of
$\Delta_{i_2}$.
\item We have $m_{i_1} = m_{i_2}$.
\endroster
\endproclaim

\secstart\label{\Charunitary} We now recall Tadi\v c's description of
unitary representations \cite{Ta} transferred to the setting of graded
Hecke algebras:

We phrase this in terms of the Kaszhdan-Lusztig classification: For any
integer $d \geq 1$ we define the $(d \times d)$-matrix
$$
n_d = \pmatrix 0 \\ 1 & 0 \\ & 1 & 0 \\ & & \ddots & \ddots \\ & & &
1 & 0 \endpmatrix.
$$
Let $L(s,e)$ be an irreducible $\HH^{1/2}_{GL_N}$-module. Then $L(s,e)$ is
unitary if and only if $(s,e)$ is conjugate to a direct sum $\bigoplus
(s_i,e_i)$ where $(s_i,e_i)$ is of one of the following forms
\roster
\item"{(I)}" $s_i = s(l_i,d_i)$ and $e_i = n_{d_i}^{\oplus l_i}$ with
$$
s(l,d) = \bigoplus_{j=1}^l\text{diag}(\frac{-l-d}2 + j,\frac{-l-d}2 + j +
1, \dots, \frac{-l+d}2 + j - 1).
$$
\item"{(II)}" $s_i = (s(l_i,d_i) + \alpha I_{l_id_i}) \oplus (s(l_i,d_i) -
\alpha I_{l_id_i})$ and $e_i = n_{d_i}^{\oplus 2l_i}$ for some real number
$\alpha_i$ with $0 < \alpha_i < \frac12$.
\endroster

\secstart\label{\unramifiedunitary} An unramified irreducible
$\HH^{1/2}_{GL_N}$-module $L(s,0)$ is unitary if and only if $s$ is
conjugated to a direct sum of diagonal matrices $s_i$ which are of one of
the following forms
\roster
\item"{(I)}" $s_i = \text{diag}(\frac{1-l_i}2,\frac{1-l_i}2 +
1,\dots,\frac{l_i-1}2)$,
\item"{(II)}" $s_i = \text{diag}(\frac{1-l_i}2 - \alpha_i, \frac{1-l_i}2 +
\alpha_i, \dots, \frac{l_i-1}2-\alpha_i,\frac{l_i-1}2 + \alpha_i)$ for some
real number $0 < \alpha_i < \frac12$.
\endroster

\secstart\label{\definesignchar} Let $V$ be a finite-dimensional Hermitian
$\HH_G$-module and let $\langle\,,\,\rangle\: V \times V \to \CC$ be a
non-degenerate Hermitian form on $V$ such that $\langle hv,v' \rangle =
\langle v,h^*v' \rangle$ for all $h \in \HH_G$ and all $v, v' \in
V$. Decompose $V = \bigoplus_{i \in I} V_i$ into an orthogonal sum of
irreducible $W$-representations. The restriction $\langle\,,\,\rangle_i$ of
$\langle\,,\,\rangle$ to $V_i$ is either positive or negative definite. For
each $\lambda \in \hat W$ we set
$$\align
\sigma^+_{\lambda}(V,\langle\,,\,\rangle) &= \#\{i \in I \mid \text{$V_i \cong
\lambda$, $\langle\,,\,\rangle_i$ is positive definite}\}, \\
\sigma^-_{\lambda}(V,\langle\,,\,\rangle) &= \#\{i \in I \mid \text{$V_i \cong
\lambda$, $\langle\,,\,\rangle_i$ is negative definite}\}, \\
\sigma_{\lambda}(V,\langle\,,\,\rangle) &=
\sigma^+_{\lambda}(V,\langle\,,\,\rangle) -
\sigma^-_{\lambda}(V,\langle\,,\,\rangle).
\endalign$$
These numbers are independent of the choice of the orthogonal decomposition
of $V$ into irreducible $W$-representations.

Assume that $V$ is irreducible as an $\HH_G$-module. In this case
$\langle\,,\,\rangle$ is uniquely determined up to a nonzero real
number. Hence the class of
$(\sigma_{\lambda}(V,\langle\,,\,\rangle))_{\lambda \in \hat W} \in
\ZZ^{\hat W}$ in $\ZZ^{\hat W}/\{\pm 1\}$ is independant of the choice of
$\langle\,,\,\rangle$ (here $\{\pm 1\}$ acts on $\ZZ^{\hat W}$ by
$\varepsilon\cdot(\sigma_{\lambda}) = (\varepsilon\sigma_{\lambda})$). We
call this class $\Sgbar(V)$.

\secstart\label{\Normalizesign} If $\Sgbar \in \ZZ^{\hat W}/\{\pm 1\}$ is the
signature character of some irreducible Hermitian representation we define
a lift $\Sigma \in \ZZ^{\hat W}$ as follows: For every irreducible
$\HH_G^{1/2}$-module $V$ there exists a unique maximal $\lambda \in \What$
such that $[V : \lambda] > 0$, and moreover we have $[V : \lambda] = 1$
\Wofstandard. We let  $\Sigma$ be the unique lift of $\Sgbar$ such that
$\Sigma_{\lambda} = 1$. Using this normalization, we get a map
$$
\Sigma\: \{\text{irreducible Hermitian $\HH_G^{1/2}$-modules}\} \arr
\ZZ^{\What}.
$$

\secstart\label{\signatureZel} The bijection $\hat W \to \hat W$ which
sends $U$ to $U \otimes \text{sgn}$ defines a $\ZZ$-linear automorphism of
$\ZZ^{\hat W}$ (by taking the corresponding permutation matrix) and this
induces a bijection of order 2 on $\ZZ^{\hat W}/\{\pm 1\}$ which we denote
again by $[(\sigma_{\lambda})] \mapsto [(\sigma_{\lambda})] \otimes
\text{sgn}$.

\proclaim{Proposition} For every irreducible $\HH_G^{1/2}$-module $V$ we have
$$
\Sgbar(\zeta(V)) = \Sgbar(V) \otimes \text{sgn}.
$$
\endproclaim

\demo{Proof} It follows directly from the definitions
that the involution ${}^*$ and the Zelevinsky involution $\zeta$ \defineZel\
commute with each other. Hence \ZelonWmodule\ implies the proposition.
\qed\enddemo


\newparagraph{Intertwining operators and the Hermitian from}

\secstart\label{\Signnotation} In the sequel we will use the
following notations: We set $G = GL_N$ and fix a standard parabolic
subgroup, given by a subset $S$ of the set of the simple roots,
corresponding to an ordered partition $(N_1,\dots,N_r)$ of $N$ and
denote by $M$ the associated standard Levi subgroup. Instead of
$\HH_S$ \Levisubalgebra\ we also write $\HH_M$.

Let $U$ be a fixed irreducible tempered representation of
$\HH^{1/2}_M$, hence it will be of the form $$ U = U_1 \boxtimes \dots
\boxtimes U_r $$ for irreducible tempered representations $U_i$ of
$\HH^{1/2}_{GL_{N_i}}$. We denote by $\langle\ ,\ \rangle_{U_i}$ the
unitary form (unique up to a positive scalar) on $U_i$ and by
$\langle\ ,\ \rangle_U$ its tensor product on $U$.

We will consider (real) Hermitian representations with Langlands data
$(S,U,\nu)$ \Langlandsclass. Hence we have for all $i$, $N_i =
N_{r+1-i}$, and $U_i$ is isomorphic to $U_{r+1-i}$. In the sequel we
choose an identification of unitary $\HH_{GL_{N_i}}$-modules $U_i \cong
U_{r+1-i}$. Finally $\nu$ will be given by an element in ${\frak
z}^{+}_S$ which we can consider as an $r$-tuple of real numbers
$\nu \in C_{(N_1,\dots,N_r)}$ with 
$$
C_{(N_1,\dots,N_r)} := \{(\nu_1,\dots,\nu_r) \in \RR^r \mid
\text{$\nu_1 > \dots > \nu_r$, $\nu_i + \nu_{r+1-i} = 0$ for $i =
1,\dots,r$}\}.
$$
We denote by $X(S,U,\nu)$ the associated standard module and by $L(S,U,\nu)$
its unique irreducible quotient, and we call $\nu$ the {\it Hermitian
parameter}.

Let $W_M$ be the Weyl group of $M$, and we set $W = W_G$. We identify $W_M$
with $S_{N_1} \times \dots \times S_{N_r}$, embedded in $W = S_N$ as
usual. Let $w_0$ be the element of maximal length in $W$ and let $w_{0,M}$
be the element of minimal length in the double coset $W_Mw_0W_M$. Note that
as $w_0$ normalizes $W_M$, we have $W_Mw_0W_M = w_0W_M = W_Mw_0$ and
$w_{0,M}^2 = 1$. Our identification of $U_i$ with $U_{r+1-i}$ then gives
an identification $U \cong w_{0,M}(U)$.

We fix the following isomorphism of $\HH_M$-modules preserving unitary
forms $$\align \tau\: U_r \boxdot \dots \boxdot U_1 &\arriso U_1
\boxdot \dots \boxdot U_r = U,\\ u_r \otimes \dots \otimes u_1 &\mapsto
u_1 \otimes \dots \otimes u_r  \endalign$$
Note that we can consider $\tau$ as an isomorphism $w_{0,M}(U) \arriso U$.

\secstart\label{\definerho} We are now introducing elements
following \cite{BM3} 1.6 and 1.7: Fix $w \in W$ and let $w =
s_1\dots s_l$ be a reduced decomposition. Then define $\rho_w =
\rho_1\rho_2\dots \rho_l$ where $\rho_i = t_{\alpha_i}\alpha_i - 2r$
if $s_i$ corresponds to the simple root $\alpha_i$. Using a result of
Lusztig ([Lu2] 5.2), it is shown in [BM3] 1.6 that $\rho_w$ does
not depend on the choice of the reduced decomposition of $w$ and that
we have for all $\theta \in \OO$
$$
\theta\rho_w = \rho_w({}^{w^{-1}}\!\theta).
\formulano\sublabel{\rhocomp}$$

\secstart\label{\definehermform} There is a unique $\OO$-linear map
$\varepsilon_M\: \HH_G \to \HH_M$ such that $\varepsilon(t_w) = t_w$
for $w \in W_M$ and $\varepsilon(t_w) = 0$ for $w \in W \setminus W_M$.

We are now going to define an Hermitian form $\beta_{S,U,\nu}$ on
$X(S,U,\nu)$ as follows: Recall \definecenter\ that for each $\mu \in
{\frak z}_S$ there exists a one-dimensional $\HH^{1/2}_M$-module,
i.e.\ a $\CC$-algebra homomorphism $\HH^{1/2}_M \to \CC$ which we
denote by $h \mapsto h(\nu)$.

For $w,w' \in W$ and $u,u' \in U \otimes \CC_{\nu}$ we set
$$
\beta_{S,U,\nu}(t_w \otimes u, t_{w'} \otimes u') =
\langle\varepsilon(t_{w^{\prime-1}}t_w\rho_{w_{0,M}})(\nu)\tau(u),u'\rangle_U.
$$

This is well defined by \rhocomp.

\secstart\label{\signatureev} The Hermitian form $(\ ,\ )$ on
$X(S,U,\nu)$ defined by
$$
(t_w \otimes u, t_{w'} \otimes u') =
\langle\varepsilon(t_{w^{\prime-1}}t_w)(\nu)u,u'\rangle_U
$$
is unitary as we have
$$
(t_w \otimes u, t_w \otimes u) = \langle u,u\rangle_U.
$$

Hence the signature character of $\beta_{\nu}$ can be computed as
follows. Consider the $\HH_G$-linear map
$$\align
A_{w_{0,M}}\: X(S,U,\nu) = \HH_G \otimes_{\HH_M} (U \otimes \CC_{\nu}) &\to
\HH_G \otimes_{\HH_M} (U \otimes \CC_{-\nu}) = X(S,U,-\nu),\\
h \otimes u &\mapsto h\rho_{w_{0,M}} \otimes \tau(u).
\endalign$$
As $\CC[W_M]$ acts trivially on $\CC_{\nu}$,
source and target of $A_{w_{0,M}}$ are canonically isomorphic as
$\CC[W]$-modules, hence we can also consider $A_{w_{0,M}}$ as an
endomorphism of a $\CC[W]$-module $X$, and for each $\lambda \in
\hat{W}$ we get an induced endomorphism $A_{w_{0,M},\lambda}$ on
$\Hom_W(\lambda,X)$. We denote the number
of positive eigenvalues minus the number of negative eigenvalues of
$A_{w_{0,M},\lambda}$ by $\sigma_{\lambda}(A_{w_{0,M}})$ and set
$\Sigma(A_{w_{0,M}}) = (\sigma_{\lambda}(A_{w_{0,M}}))_{\lambda \in
\What}$. Then we have
$$
\Sigma(A_{w_{0,M}}) = \Sigma(\beta_{\nu}).
$$

\proclaim{Proposition} Assume that $A_{w_{0,M}}$ is nonzero. Then its image
is isomorphic to $L(S,U,\nu)$ and $\beta_{\nu}$ is up to a scalar the
Hermitian form given by the involution ${}^*$ \hermitianinvol.
\endproclaim

\demo{Proof}: Let $X' \subset X(S,U,\nu)$ be the maximal submodule. Then we
ave $X/X' = L(S,U,\nu)$ \subofstandard. We have to show that the
restriction $A'\: X' \arr X(S,U,\nu)^h$ of $A_{w_{0,M}}$ to $X'$ is always
zero. As $X(S,U,\nu)$ has a unique irreducible quotient, namely
$L(S,U,\nu)$, $X(S,U,\nu)^h$ has a unique irreducible submodule, namely
$L(S,U,\nu)^h$ which is isomorphic to $L(S,U,\nu)$ as $L(S,U,\nu)$ is
Hermitian. As $X(S,U,\nu)^h$ is of finite length, the image of $A'$ has to
contain this unique irreducible submodule if $A'$ is nonzero. This would
imply that $X'$ has a subquotient which is isomorphic to $L(S,U,\nu)$ but
this is a contradiction to the fact that $L(S,U,\nu)$ occurs with
multiplicity one in $X(S,U,\nu)$.
\qed\enddemo


\break

\newparagraph{Reducibility walls}

\secstart\label{\Sigmalocconst} We keep the notations from \Signnotation.

For fixed $S$ und $U$ we call $\nu \in C_{(N_1,\dots,N_r)}$ {\it
irreducible} if the $\HH^{1/2}_{GL_N}$-module $X(S,U,\nu)$ is irreducible and
denote by $C_{(N_1,\dots,N_r)}^0$ the set of irreducible $\nu$.

Every irreducible Hermitian irreducible $\HH^{1/2}_{GL_N}$-module $V$
defines a signature character $\Sigma(V) \in \ZZ^{\hat{W}}$ \Normalizesign. We
obtain a map
$$\align
C_{(N_1,\dots,N_r)} &\to \ZZ^{\hat{W}},\\
\nu &\mapsto \Sigma(L(S,U,\nu)).
\endalign$$

By \definehermform\ and \signatureev\ this map is locally constant on
$C_{(N_1,\dots,N_r)}^0$.

\secstart Write $L(S,U,\nu) = L(\Delta_1,\dots,\Delta_s)$ \BZclass. Note
that the number and the length of the segments do not depend on
$\nu$. The standard module $X(S,U,\nu)$ is irreducible if and only if no
two of the segments can be linked. Hence the reducibility locus,
i.e. $C_{(N_1,\dots,N_r)} \setminus C_{(N_1,\dots,N_r)}^0$, is a union of
hyperplanes of the form
$$
H_{\alpha} = \{\nu \in C_{(N_1,\dots,N_r)} \mid \langle \alpha, \nu +
\chi \rangle = 1\}
$$
where $\alpha \in R^+$ runs through a certain set of positive roots (cf.\
(9.4)) and where $\chi$ is the central character of the tempered
representation $U$.

\secstart We make this more concrete in the case of an unramified
representations, i.e.\ for the case $S = \emptyset$ and hence $U$ the
trivial representation. In this case all segments have length $1$ and
$\nu = (\nu_1,\dots,\nu_N)$ is irreducible if and only if $\langle
\nu, \alpha \rangle \not=1$ for all $\alpha \in R^+$, i.e.~$\nu_i -
\nu_j \not= 1$ for all $i < j$.

For every root $\alpha \in R^+$ we define the corresponding
reducibility wall $$ H_{\alpha} = \{\nu \in C_N \mid \langle \alpha,
\nu \rangle = 1\}.  $$

Set $M = [N/2]$. Via $(\nu_1,\dots,\nu_N) \mapsto (\nu_1,\dots,\nu_M)$
we can identify $C_N$ with $D_M$ where $$ D_M = \{(x_1,\dots,x_M) \in
\RR^M \mid x_1 > \dots > x_M > 0\}.  $$

Via this identification the reducibility walls can be described as
follows. Assume first that $N$ is even. Here we have that all nonempty
reducibility walls are the following: $$\alignat2 H_{ij}^\pm &=
\{(x_1,\cdots,x_M)\in D_M \mid x_i \pm x_j = 1\},&&\qquad \text{for
$1\leq i<j\leq M$} \\ \intertext{and} H_{i,\frac 1 2} &=
\{(x_1,\cdots,x_M)\in D_M \mid x_i =  \frac 1 2\}, &&\qquad\text{for
$1\leq i \leq M$}.\\ \intertext{For $N$ odd, we have in addition to
the walls above those of the form} H_{j,1} &= \{(x_1,\cdots,x_M)\in
D_M \mid x_j = \pm 1\},&&\qquad\text{for $1\leq j \leq M$}.
\endalignat$$

We will give a description of the reducibility walls in terms of the
roots, namely, we will determine a subset of the set of roots such
that each reducibility plane is determined by one element in this
subset. The first easy remark is that $$ \aligned H^+_\alpha &=
H_{-\alpha}^-\\ H^-_\alpha &= H_{-\alpha}^+ \endaligned $$ for any
$\alpha$ in the root lattice. Denote by  $$\aligned K =
\{\alpha_{i,j}^\pm = e_i\pm e_j \,| & \, 1\leq i<j\leq N\}\cup
\{\alpha_i=e_i-e_{N-i+1}\,|\, i=1,\cdots,M\}\\
 & \cup \{\alpha_{i,M}=e_i-e_{M+1} \,|\,  i=1,\cdots,M\}.  \endaligned
$$

We have a one to one correspondance between the set of all vanishing
walls and elements in $K$, namely  $H_{i,j}^\pm$ maps to
$\alpha_{i,j}^\pm$, and also between each $H_{i,\frac 1 2}$ to
$\alpha_i=e_i-e_{N-i+1}$, if $N$ is even. Now if $N$ is odd, we also
have that $H_{j,1}^\pm$ corresponds to $\alpha_{i,M}=e_i-e_{M+1}$.

We have the following trivial intersection rules for the walls:
\roster
\item $H_{ij}^- \cap H_{ik}^- = H^-_{ik} \cap H^-_{jk} =
H_{ij}^+ \cap H_{ik}^+ = H^+_{ik} \cap H^+_{jk} = \emptyset$.
\item $H_{ij}^+ \cap H_{jk}^+ = \emptyset$.
\item $H^+_{ij} \cap H^-_{ij} = \emptyset$.
\item $H_{i,\frac12} \cap H_{j,\frac12} = H_{i,1} \cap H_{j,1}= \emptyset$.
\item $H^+_{ij} \cap H^-_{jk} = \emptyset$.
\item $H^-_{ij} \cap H_{i,\frac12} = H^-_{ij} \cap H_{i,1} = \emptyset$.
\item $H^+_{ij} \cap H_{i,\frac12} = H^+_{ij} \cap H_{j,\frac12} = H^+_{ij}
\cap H_{i,1} = H^+_{ij} \cap H_{j,1} = \emptyset$.  \endroster

Now let $R' \subset R^+$ be a set of positive roots. We set $$ H_{R'}
= \bigcap_{\alpha \in R'}H_{\alpha}.  $$

\proclaim{Proposition} Let $R' \subset R^+$ be such that $H_{R'}$
consists of a single point $\nu = (\nu_1,\dots,\nu_N)$. Then $L(\nu)$
is unitary.
\endproclaim

\demo{Proof} Let $J \subset R'$ be a minimal subset such that $H_J =
H_{R'}$. In particular we have $\#J = M$. We say that two walls $H$ and
$H'$ in $J$ are equivalent if there exist walls $H = H_0,
H_1,\dots,H_{r-1}, H_r = H'$ in $J$ such that the set of indices for
the wall $H_k$ has nonempty intersection with the set of indices of
$H_{k+1}$. Let $J = \bigcup J_{\alpha}$ be the decomposition into
equivalence classes with respect to this relation. $J_{\alpha}$ has to be a
set of walls of one of the following forms:
\roster
\item"{(i)}" $J_{\alpha} =
\{H^-_{i_1,i_2}, H^-_{i_2,i_3}, \dots, H^-_{i_{m-1},i_m}\}$.
\item"{(ii)}" $J_{\alpha} = \{H^-_{i_1,i_2}, H^-_{i_2,i_3}, \dots,
H^-_{i_{m-1},i_m},H_{i_m,\frac12}\}$.
\item"{(iii)}" $J_{\alpha} =
\{H^-_{i_1,i_2}, H^-_{i_2,i_3}, \dots, H^-_{i_{m-1},i_m},H_{i_m,1}\}$.
\item"{(iv)}" $J_{\alpha} = \{H^-_{i_1,i_2}, H^-_{i_2,i_3}, \dots,
H^-_{i_{m-1},i_m},H^+_{i_m,i_{m+1}}\}$.
\endroster

As we have $\# J = \sum \# J_{\alpha} = M = \#\{\text{indices of walls
occuring in $J$}\}$, in $J_{\alpha}$ have to be $\#J_{\alpha}$
indices, and this means that only the cases (ii) and (iii) are
possible. As the $\nu_i$ are pairwise different, each of the cases
(ii) and (iii) can occur only once. Now it follows from Tadi\v c's
description of unitary representations \unramifiedunitary\ that
$L(\emptyset,{\bold 1}, \nu)$ has to be unitary.
\qed\enddemo


\newparagraph{Signature at Infinity}

\secstart Recall from \Sigmalocconst\ that the signature character is
constant on the connected components of $C_{(N_1,\dots,N_r)}^0$. It is
easy to see that there is a unique connected component $C_{\infty}$ of
$C^0_{(N_1,\dots,N_r)}$ such that for all roots $\alpha \in \Pi
\setminus S$ the set of real numbers $\{\nu(\alpha) \mid \nu \in C\}$
is not bounded. In this section we want to describe the signature of
$X(S,U,\nu) = L(S,U,\nu)$ for $\nu \in C_{\infty}$. For this we
follow Barbasch and Moy \cite{BM3}.

\secstart\label{\definesigninf} As a $W$-module we have $X(S,U,\nu)
\cong \CC[W] \otimes_{\CC[W_M]} U$. Let $\varepsilon\: \CC[W] \to
\CC[W_M]$ be the $\CC$-linear map with $\varepsilon(t_w) = t_w$ for $w
\in W_M$ and $\varepsilon(t_w) = 0$ for $w \in W \setminus W_M$ and
define a $\CC[W]$-invariant Hermitian form on $X(S,U,\nu)$ by
$$
\beta_{\infty}(t_w \otimes u, t_{w'} \otimes u') = \langle
\varepsilon(t_{w^{\prime-1}}t_wt_{w_{0,M}})\tau(u),u'\rangle_U.
$$
As in \definesignchar\ this defines a signature character
$\Sigma_{\infty} \in \ZZ^{\hat{W}}$ which depends only on $S$
und $U$.

A limit argument using \Sigmalocconst\ (cf.\ \cite{BM3} 2.3)
shows that for $\nu \in C_{\infty}$ the signature character of
$\beta_{\nu}$ and the signature character of $\beta_{\infty}$
coincide. Hence we get that for $\nu \in C_{\infty}$ the class of
$\Sigma_{\infty}$ in $\ZZ^{\What}/\{\pm 1\}$ is equal to the class of the
signature of $\beta_{\nu}$ in $\ZZ^{\What}/\{\pm 1\}$.

\secstart\label{\signbyendo} To calculate the signature character of
$\beta_{\infty}$ we make the following general remark. Let
$(V,\langle\,,\,\rangle)$ be any finite-dimensional complex unitary space
and $f \in \End(V)$ be a self adjoint endomorphism (in particular $f$ is
semisimple). We define a Hermitian form $\beta_f$ on $V$ by
$$\beta_f(v,v') = \langle fv, v' \rangle.$$
This form is non-degenerate iff $f$ is invertible and in this case its
signature is equal to
$$
\#\{\text{positive eigenvalues of $f$}\} -
\#\{\text{negative eigenvalues of $f$}\}
$$
where we count eigenvalues with multiplicity. Now assume further that $f^2
= \text{id}_V$. Then $f$ has only eigenvalues $+1$ or $-1$, and the
signature of $\beta_f$ is nothing but $\Tr(f)$.

\secstart\label{\calculateinfinity} Let us apply this to the unitary form
$\langle t_w \otimes u,t_{w'} \otimes u'\rangle =
\langle\varepsilon(t_{w^{\prime-1}}t_w)u,u'\rangle_U$ on $X(S,U,\nu)$ and
to $f = r_{w_{0,M}}$ given by
$$
r_{w_{0,M}}\: t_w \otimes u \mapsto t_wt_\w0M \otimes \tau(u).
$$
Then it follows that $\Sigma_{\infty} = (\sigma_{\infty,\lambda}) \in
\ZZ^{\hat{W}}$ with
$$
\sigma_{\infty,\lambda} =
\dim(\lambda)^{-1}\Tr(r_{w_{0,M}}|_{X_{\lambda}})\formulano
\sublabel{\signwithtrace}
$$
where $\dim(\lambda)$ denotes the complex dimension of the irreducible
$W$-representation corresponding to $\lambda \in \hat{W}$ and $X_{\lambda}$
denotes the $\lambda$-isotypical component of the left $W$-module $X =
X(S,U,\nu)$.

Now let $\chi_{\lambda}$ be the character on $W$ corresponding
to the irreducible representation $\lambda$. The projection $p_{\lambda}$
from $X$ onto its isotypical component $X_{\lambda}$ is given by
$$
t \otimes u \mapsto \frac{\dim(\lambda)}{\#W}\sum_{w\in W}\chi_{\lambda}(w)(wt
\otimes u).
$$
Hence if we define for $w \in W$ a $\CC$-linear endomorphism $f^w$ of $X$ by
$$t \otimes u \mapsto wt\w0M \otimes \tau(u),$$
we have
$$\align
\sigma_{\infty,\lambda} &= \dim(\lambda)^{-1}\Tr(p_{\lambda} \circ r_\w0M)\\
&= (\#W)^{-1}\sum_{w \in W}\chi_{\lambda}(w)\Tr(f^w).
\endalign$$

For $z \in W$ let $\ell_z\: X \arr X$ be the left multiplication with
$z$. Then we have
$$
\ell_z \circ f^w \circ \ell_z^{-1} = f^{zwz^{-1}}.
$$
In particular we see that $\Tr(f^w)$ depends only on the conjugacy class of
$w$. Identifying conjugacy classes in $W$ with irreducible representations of
$W$, we get
$$
\sigma_{\infty,\lambda} = (\#W)^{-1}\sum_{\mu \in
\hat{W}}N(\mu)\chi_{\lambda}(\mu)\Tr(f^{w(\mu)})\formulano
\sublabel{\signexpression}
$$
where $N(\mu)$ is the number of elements of $W$ in the conjugacy class
$\mu$ and where $w(\mu)$ is some element in $\mu$.

\secstart It remains to calculate the trace of $f^w$. Let $[W/W_M]$ be a
system of representatives in $W$ of the quotient $W/W_M$. As a $\CC$-vector
space, $X$ is isomorphic to $\bigoplus_{w \in [W/W_M]}\CC t_w \otimes
U$. Hence
$$
\Tr(f^w) = \sum_{z \in [W/W_M], \atop z^{-1}wz\w0M \in W_M}
\Tr(z^{-1}wz\w0M\tau|U).\formulano\sublabel{\traceexpression}
$$

It remains to determine for $x = (x_1,\dots,x_r) \in W_M = S_{N_1} \times
\dots S_{N_r}$ the trace of the endomorphism
$u \mapsto x\tau(u)$ of $U$. For this we set $U'_i = U_i \otimes
U_{r+1-i}$ for $i = 1,\dots,[r/2]$ and also $U'_{(r+1)/2} = U_{(r+1)/2}$ if
$r$ is odd. Then we have
$$
U = U'_1 \otimes \dots \otimes U'_{[(r+1)/2]}
$$
and the endomorphism $x\tau\vert U$ is the tensor product of the
endomorphisms $x'_i\tau_i\vert U'_i$ where $x'_i = (x_i,x_{r+1-i})$ (and
$x'_{[(r+1)/2]} = x'_{[(r+1)/2]}$ if $r$ is odd) and $\tau_i$ is the
endomorphism of $U'_i$ which switches the two components (and is the
identity on $U'_{(r+1)/2}$ if $r$ is odd). Therefore we have
$$
\Tr(x\tau\vert U) = \prod_{i=1}^{[(r+1)/2]}\Tr(x'_i\tau_i\vert
U'_i).\formulano\sublabel{\decomposetrace}
$$

\secstart Hence we are reduced to the following situation: For an integer
$M \geq 1$ let $V$ be a finite-dimensional  $S_M$-module und set $U = V
\otimes V$. For $(x_1,x_2) \in S_M \times S_M$ we want to determine the
trace of $\ell_{(x_1,x_2)} \circ \tau$ with
$$\align
\tau\: v_1 \otimes v_2 &\mapsto v_2 \otimes v_1,\\
\ell_{(x_1,x_2)}\: v_1 \otimes v_2 &\mapsto x_1v_1 \otimes x_2v_2.
\endalign$$

We claim that
$$
\Tr(\ell_{(x_1,x_2)} \circ \tau) = \sum_{\lambda \in \hat{S}_M}[V :
\lambda]\chi_{\lambda}(x_1x_2).
$$

First note that if we have a decomposition of $S_M$-modules $V = V' \oplus
V''$, a concrete matrix multiplication shows that
$$
\Tr((\ell_{(x_1,x_2)} \circ \tau)|(V \otimes V))
= \Tr((\ell_{(x_1,x_2)} \circ \tau)|(V' \otimes V'))
  + \Tr((\ell_{(x_1,x_2)} \circ \tau)|(V'' \otimes V'')).
$$
Therefore we can assume that the $S_M$-module $V$ is irreducible, say of
isomorphism class $\lambda \in \hat{S}_M$. Of course, then we know that the
trace of left multiplication with $x_1x_2$ on $V$ is simply
$\chi_{\lambda}(x_1x_2)$. Hence the claim follows if we prove the following
lemma:

\proclaim{Lemma} Let $k$ be a field, let $V$ be a finite-dimensional
$k$-vector space and let $f_1$ and $f_2$ be two endomorphisms of $V$ and let
$\tau\: V \otimes V \to V \otimes V$ the map that switches components. Then
we have
$$
\Tr((f_1 \otimes f_2) \circ \tau) = \Tr(f_1 \circ f_2).
$$\endproclaim

\demo{Proof} One easily reduces to the case that $k$ is algebraically
closed. Using the fact that the set of semisimple endomorphisms is dense
within the sapce all endomorphisms with respect to the Zariski topology, we
can also assume that $f_1$ is semisimple. Choose a basis $(e_i)$ of $V$ such
that $f_1$ is given by a diagonal matrix $A_1$ with respect to this
basis. Denote by $A_2$ the matrix of $f_2$. Then the map $(f \otimes f') \circ
\tau$ sends $e_i \otimes e_j$ to $(A_1)_{jj}e_j \otimes
(\sum_{l}(A_2)_{li}e_l)$. Hence we see that the trace of $(f_1 \otimes f_2)
\circ \tau$ is nothing but $\sum_{i}(A_1)_{ii}(A_2)_{ii}$ which is the same as
the trace of $f_1 \circ f_2$.
\qed\enddemo

\secstart\label{\unramifiedinfinity} In \calculateinfinity\ we described
an algorithm which reduced the computation of the signature character of
$X(S,U,\nu)$ for $\nu \in C_{\infty}$ to the computation of the isomorphism
class of $U$ as a $\CC[W_M]$-module (which we know by \standWmodule) and
the computation of the characters of the symmetric group (which is
well-known).

In the case of unramified representations (i.e.~where $S$ is empty and $U$
is the trivial representation) we have:

\proclaim{Corollary} For unramified representations the signature character of
$\beta_{\infty}$ is equal to $(\sigma_{\infty,\lambda})_{\lambda \in
\hat{W}} \in \ZZ^{\hat{W}}$ with
$$
\sigma_{\infty,\lambda} = \chi_{\lambda}(w_0).
$$
\endproclaim

\demo{Proof} This is a straight forward application of the algorithm
above. Note that we do not only have equality in $\ZZ^{\hat{W}}/\{\pm 1\}$
but even in $\ZZ^{\hat{W}}$ as we have $\chi_{\bold 1}(w_0) = 1$.
\qed\enddemo


\newparagraph{Wall crossing}

\secstart We continue to use the notations of \Signnotation. Also denote by
$\chi$ the central character of $U$. Let $(s,e)$ be a pair as in \KLclass\
corresponding to $L(S,U,\nu)$. We can conjugate $(s,e)$ such that $s = \chi
+ \nu \in \frak{t}$. Note that the conjugacy class of $e$ depends only on
$S$ und $U$. Let $R_M \subset R$ be the root system of the Levi subroup
$M$. We further assume that we are not in the trivial case that the cone
$C_{(N_1,\dots,N_r)}$ of Hermitian parameters is empty.

The tempered representations $U_i$ are of the form
$L(\Delta^i_1,\dots,\Delta^i_{r_i})$ where $\Delta^i_j = \Delta(-(m^i_j -
1)/2,m^i_j)$ is a segment with center 0. We order these segments such that
$j \mapsto m^i_j$ is nonincreasing.

Let $\Delta_1,\dots,\Delta_m$ be a tuple of segments such that
$L(\Delta_1,\dots,\Delta_m) = L(S,U,\nu)$. We call the number of pairs
$(\Delta_i,\Delta_j)$ such that $\Delta_i$ precedes $\Delta_j$ the {\it
height} of $\nu$. Hence $\nu$ is irreducible if and only if the height of
$\nu$ is zero. For every reducibility wall $H$ the height of $\nu \in H$ is
constant equal to some positive integer $\text{height}(H)$ outside those
points that lie on an intersection of $H$ with some other reducibility
wall. We call $\text{height}(H)$ the {\it height} of the reducibility wall.

Our goal in this chapter is to cross these walls. More precisely let
$\nu_0$ be a Hermitian parameter, such that there exists a unique
reducibility wall $H$ with $\nu_0 \in H$. Then we can find a small positive
real number $\varepsilon$ and a non constant map
$$
(-\varepsilon, \varepsilon) \to C_{(N_1,\dots,N_r)},\quad t \mapsto \nu(t)
$$
which is the restriction of an affine map $\RR \to \RR^N$ such that
\roster
\item"{(i)}" $\nu(0) = \nu_0$.
\item"{(ii)}" $\nu(t)$ is irreducible for all $t \not= 0$.
\endroster

Then for $t > 0$ (resp.\ $t < 0$), the signature character of
$L(S,U,\nu(t))$ is constant and we call it $\Sigma^+$ (resp.\
$\Sigma^-$). The goal of this chapter is to find an expression for
$\Sigma^+ + \Sigma^-$ and for $\Sigma^- - \Sigma^+$.

\secstart\label{\Jantzenfiltration} To do this we use a filtration which
is called in [Vo] 3 the Jantzen filtration: Let $\HH$ be a $\CC$-algebra with
a complex antilinear involution ${}^*$ and let $E$ be a left $H$-module
which is a finite dimensional $\CC$-vector space. Let $\beta_t$ be a real
analytic family of Hermitian forms on $E$ defined for small real $t$ such
that $\beta_t(hx,y) = \beta_t(x,h^*y)$ for all $h \in \HH$, $x,y \in
E$. Assume that $\beta_t$ is non-degenerate for $t \not= 0$. Then there is
a unique sequence of subspaces
$$E = E^0 \supset E^1 \supset \dots \supset E^n = (0)$$
such that the meromorphic familiy of Hermitian forms $\beta^i_t=
\frac{1}{t^i}\beta_t(x,y)\vert_{E^i}$ can be extended to an analytic family
of Hermitian forms on $E^i$ and such that the radical of $\beta^i_0$ is
equal to $E^{i+1}$. In particular $\beta^i_0$ induces a non-degenerate
pairing $\beta^i$ on $\gr^i(E) = E^i/E^{i+1}$.

Explicitly $E^i$ can be defined as
$$E^i = \{x \in E \mid \text{$\beta_t(x,y)$ vanishes at least to order $i$
at $t = 0$, for any $y \in E$}\}.$$

Note that our definition of the $E^i$ varies a little bit from the
definition given in loc.\ cit., but it is easily seen that both definitions are
equivalent. Moreover it follows from the definition that the $E^i$ are also
$\HH$-submodules of $E$ and that $\beta^i_t$ are Hermitian with respect to
the involution ${}^*$.

Let $\sigma^i$ be the signature of $\beta^i$, and denote by $\sigma^+$
(resp.\ $\sigma^-$) the signature of $\beta_t$ for small positive (resp.\
negative) $t$. Then we have (loc. cit.)
$$\align
\sigma^+ - \sigma^- &= 2\sum_{\text{$i \geq 1$ odd}}\sigma^i,\\
\sigma^+ + \sigma^- &= 2\sum_{\text{$i \geq 0$ even}}\sigma^i.\\
\endalign$$

\secstart\label{\orddet} Now let $\langle\ ,\ \rangle$ be a fixed unitary
form on $E$ and let $f_t\: E \to E$ be the analytic family of selfadjoint
endomorphisms such that $\beta_t(x,y) = \langle f_t(x), y\rangle$. Then
$\det(f_t)$ is non-zero for $t \not= 0$ and the analytic map $t \mapsto
\det(f_t)$ has a zero of order
$$
\sum_{i > 0} i(\dim(E^i) - \dim(E^{i-1})) = \sum_{i > 0} i\dim(\gr^i(E))
$$
in $t = 0$.

\secstart\label{\distinguishsides} Let $\bar{\Cal{O}_e}$ be the closure of
the orbit of $e$ in $\Cal{N}(s) = \{n \in \frak{g} \mid [s,n] = n\}$ under
the action of $Z_G(s)$. For all $\nu$ we define $R_{\nu}$ as the set of
roots $\alpha \in R^+$ such that $\frak{g}^{\alpha} \subset \Cal{N}(s)$ and
such that there exists a nonzero element $e_{\alpha} \in \frak{g}^{\alpha}$
such that $e + e_{\alpha} \notin \bar{\Cal{O}_e}$. These are those roots
which ``link'' segments of $X(S,U,\nu)$. We have $R_{\nu} = \emptyset$ if
and only if $\nu$ is irreducible. If $\nu$ lies on a unique reducibilty
wall $H$, $R_{\nu}$ is independent of $\nu$, and we call it $R_H$. 

Then for all $\alpha \in R_H$, $\sgn(\alpha(\nu(t)+\chi)-1)$ is independent of
$\alpha \in R_H$ and we have $\sgn(\alpha(\nu(t)+\chi)-1) = -
\sgn(\alpha(\nu(-t)+\chi)-1)$ for all $t \in (-\varepsilon,
\varepsilon)$. After a possible substitution of $t \mapsto -t$ we can and
will from now on assume that $\sgn(\alpha(\nu(t)+\chi)-1)) = \sgn(t)$.

\secstart\label{\Heightone} Let $H$ be a reducibility wall of height
one and fix $t \mapsto \nu(t)$ as above. We write $L(S,U,\nu(0))$ as
$L(\frak{M})$ where $\frak{M}$ is a multiset of segments
$\Delta_1,\dots,\Delta_m$. As $H$ is of height one, there exists a unique
pair $(\Delta_i,\Delta_j)$ of segments, such that $\Delta_i$ precedes
$\Delta_j$. We denote $\frak{M}'$ the multiset consisting of the segments
$\Delta_l$ for $l \not= i,j$ and the segments $\Delta_i \cap \Delta_j$ and
$\Delta_i \cup \Delta_j$. By \BZclass\ and \subofstandard\ we have an exact
sequence
$$
0 \to L(\frak{M}') \to X(\frak{M}) \to L(\frak{M}) \to 0
\formulano\sublabel{\exactone}
$$
and $L(\frak{M}')$ is again Hermitian, and we have $L(\frak{M}') =
X(\frak{M}')$ otherwise $\nu(0)$ would not lie on a unique reducibility wall.

We now apply \Jantzenfiltration\ to $E = X(\nu(t))$ which as a
$\CC[W]$-module does not depend on $t$ and its Hermitian form $\beta_t$. 
It is of the form
$$
X(\nu(t)) = E^0 \varsupsetneqq L(\frak{M}') = E^1 = \dots = E^{\omega}
\varsupsetneqq (0) = E^{\omega+1}
$$
for some $\omega \geq 1$, and we have
$$
\gr^0(E) = L(\nu(0)).
$$
As $L(\frak{M}')$ is an irreducible module, its signature
character has to be equal to the signature character of $\beta_{\omega}$ up
to a sign $\varepsilon$.

If we set $\Sigma^+ =  \lim_{t \to 0^+} \Sigma(L(S,U,\nu(t)))$ and
$\Sigma^- = \lim_{t \to 0^-} \Sigma(L(S,U,\nu(t)))$ we will therefore have
$$\align
\Sigma^+ + \Sigma^-
&= \cases
2\Sigma(L(S,U,\nu(0))),&\text{if $\omega$ is odd,}\\
2\Sigma(L(S,U,\nu(0))) + 2\varepsilon\Sigma(L(\frak{M}')),&\text{if
$\omega$ is even,}
\endcases\\
\Sigma^+ - \Sigma^-
&= \cases
2\varepsilon\Sigma(L(\frak{M}')),&\text{if $\omega$ is odd,}\\ 0,&\text{if
$\omega$ is even.}
\endcases
\endalign$$

Hence if $\omega$ is even, the signature characters on both sides of the
wall are equal. If $\omega$ is odd, we can calculate $\Sigma^+$ (resp.\
$\Sigma^-$) if we know $\Sigma(L(\frak{M}'))$ and $\Sigma^-$ (resp.\
$\Sigma^+$) as this allows us also to calculate $\varepsilon$: By
\standWmodule\ we know that the sign representation occurs with
multiplicity one in $X(S,U,\nu(t))$ for all $t$ and $L(\frak{M}')$. Hence
the equality
$$
(\Sigma^+)_{\sgn} - (\Sigma^-)_{\sign} =
2\varepsilon\Sigma(L(\frak{M}'))_{\sgn}
$$
implies that we have $(\Sigma^+)_{\sgn} = -(\Sigma^-)_{\sign}$ and
$$
\varepsilon = \Sigma(L(\frak{M}'))_{\sgn}(\Sigma^+)_{\sgn}.
$$

\secstart\label{\Conjectureone} In fact we conjecture the following:

\proclaim{Conjecture for reducibility walls of height one} We always have
the $\omega$ is odd. Hence
$$\align
\Sigma^+ - \Sigma^- &= \varepsilon2\Sigma(L(\frak{M}')),\\
\Sigma^+ + \Sigma^- &= 2\Sigma(L(\frak{M})).
\endalign$$
with $\varepsilon = \Sigma(L(\frak{M}'))_{\sgn}(\Sigma^+)_{\sgn} = -
\Sigma(L(\frak{M}'))_{\sgn}(\Sigma^-)_{\sgn}$.

\endproclaim

\secstart\label{\Heighttwo} Now let $H$ be a reducibility wall of
height two and again fix $t \mapsto \nu(t)$ as above. We write
$L(S,U,\nu(0))$ as $L(\frak{M})$ where $\frak{M}$ is a multiset of segments
$\Delta_1,\dots,\Delta_m$. Now there exist two pairs
$(\Delta_{i_1},\Delta_{j_1})$ and $(\Delta_{i_2},\Delta_{j_2})$ of
segments, such that $\Delta_{i_k}$ precedes $\Delta_{j_k}$ for $k =
1,2$. We can assume that $\Delta_{j_2}$ does not precede
$\Delta_{i_1}$. Note that it is possible that $\Delta_{j_1}$ precedes
$\Delta_{i_2}$. We denote $\frak{M}_1$ (resp.\ $\frak{M}_2$, resp.\
$\frak{M}'$) the multiset of segments which we get from $\frak{M}$ by
linking $(\Delta_{i_1}, \Delta_{j_1})$ (resp.\
$(\Delta_{i_2},\Delta_{j_2})$, resp.\ both $(\Delta_{i_1}, \Delta_{j_1})$
and $(\Delta_{i_2}, \Delta_{j_2})$.

By \subofstandard\ we know that $L(\frak{M})$ is the unique irreducible
quotient of $X(\frak{M})$, $L(\frak{M'})$ is the unique irreducible
submodule of $X(\frak{M})$ and that the other irreducible subquotients of
$X(\frak{M})$ are isomorphic to $L(\frak{M}_1)$ and
$L(\frak{M}_2)$. We assume that $L(\frak{M}_1)$ and
$L(\frak{M}_2)$ both occur with multiplicity one $X(\frak{M})$.
The standard modules satisfy the inclusions $X(\frak{M}') \subset
X(\frak{M}_i) \subset X(\frak{M})$ for $i = 1,2$. $L(\frak{M}') =
X(\frak{M}')$ is the unique irreducible submodule of $X(\frak{M})$, we have
exact sequences
$$
0 \to L(\frak{M}') \to X(\frak{M}_i) \to L(\frak{M}_i) \to 0,
$$
and $L(\frak{M}) = X(\frak{M})/(X(\frak{M}_1) + X(\frak{M}_2))$.

Further $L(\frak{M}')$ is Hermitian and as $\nu(0)$ lies only on a unique
reducibility wall, we have $L(\frak{M}') = X(\frak{M}')$. $L(\frak{M}_1)$
is non-Hermitian, its Hermitian dual is $L(\frak{M}_2)$.

Again we apply \Jantzenfiltration\ to $E = X(\frak{M})$ with its family of
Hermitian forms $\beta_t$. As the Hermitian forms induced on the graded
pieces of the Jantzen filtration are non-degenerate, it follows that the
Jantzen filtration will be of the form
$$\split
X(\frak{M}) = E^0 &\varsupsetneqq X(\frak{M}_1) + X(\frak{M}_2) = E^1
= \dots = E^{\omega_1} \\
&\varsupsetneqq X(\frak{M}') = E^{\omega_1+1} = \dots
= E^{\omega_2} \\
&\varsupsetneqq (0) = E^{\omega_2+1}
\endsplit$$
for integers $\omega_2 \geq \omega_1 > 0$. Here $\omega_2 = \omega_1$ means
that only $L(\frak{M})$ and $X(\frak{M}_1) + X(\frak{M}_2)$ occur as graded
pieces in the filtration.

If $\omega_2 > \omega_1$, the signature character of the form
$\beta^{\omega_2}$ induced on the irreducible module $\gr^{\omega_2}(E) =
L(\frak{M}')$ is equal to the signature character of $L(\frak{M}')$ up to a
sign $\varepsilon$.

Again we can use \Jantzenfiltration\ to compute difference and sum of
$\lim_{t\to 0^+}\Sigma(L(S,U,\nu(t)))$ and $\lim_{t\to
0^-}\Sigma(L(S,U,\nu(t)))$ in terms of $\Sigma(L(\frak{M}))$,
$\Sigma(L(\frak{M}'))$ and $\Sigma(L(\frak{M}_1) \oplus
L(\frak{M}_2))$. Note that the last signature is zero as $L(\frak{M}_1)$ and
$L(\frak{M}_2)$ are dual to each other.

\secstart\label{\Prepareproof} We now consider wall crossing in the
unramified case. Hence from now on we assume that $S$ is empty. Hermitian
representations are then given by elements $\nu \in C_{(1,\dots,1)} =
\{(\nu_i) \in \RR^N \mid \nu_i + \nu_{N+1-i} = 0\}$. For simplicity we
write $X(\nu)$ for the corresponding standard module and $L(\nu)$ for its
unique irreducible quotient.

We now consider the $\HH_G$-linear homomorphism
$$
A_{w_0}\: h \otimes 1 \mapsto h\rho_{w_0} \otimes 1
$$ 
defined in \signatureev. As the leading term of $\rho_{w_0}$ is
$t_{w_0}\prod_{\alpha \in R^+}\alpha$, $A_{w_0}$ is nonzero as we have
$\prod_{\alpha \in R^+}\langle \alpha,\nu \rangle \not= 0$. Hence we can
calculate the signature character of $L(\nu)$ using the form $\beta_{\nu}$
defined in \definehermform. If we write $w_0 = s_ks_{k-1}\dots s_1$ as a
product of simple reflections, we get a decomposition $A_{w_0} = T_k \circ
T_{k-1} \circ \dots \circ T_1$ where $T_i$ is an $\HH_G$-linear map
$$
T_i\: \HH_G \otimes_{\HH_T} \CC_{s_{i-1}\dots s_1(\nu)} \arr \HH_G
\otimes_{\HH_T} \CC_{s_i\dots s_1(\nu)}.
$$
Source and target of each $T_i$ are canonically isomorphic as
$\CC[W]$-modules and hence we can consider $T_i$ as an endomorphism of
$\CC[W]$-modules.

For each simple root $\alpha \in \Pi$ we define the Levi subgroup
$G^{\alpha}$ of $G$ by $\Lie(G^{\alpha}) = \frak{t} \oplus
\frak{g}^{\alpha} \oplus \frak{g}^{-\alpha}$. If the simple reflection
$s_i$ corresponds to the simple root $\alpha_i$, $T_i$ can be written as
$\id_{\HH_G} \otimes T_i^{\alpha_i}$ with
$$\align
T_i^{\alpha_i}\: \HH_{G^{\alpha_i}} \otimes_{\HH_T} \CC_{s_{i-1}\dots
s_1(\nu)} &\arr \HH_{G^{\alpha_i}} \otimes_{\HH_T} \CC_{s_i\dots s_1(\nu)},\\
h \otimes 1 &\mapsto h(t_{\alpha_i}\alpha_i - 1) \otimes 1.
\endalign$$
Similarly as above we can consider $T_i^{\alpha_i}$ as an endomorphism of
$\CC[S_2]$-modules. Further by \standWmodule\ we know source and target are
isomorphic to $\CC[S_2] = \sign \oplus \bold{1}$. An easy calculation shows
that $T_i^{\alpha_i}$ acts on $\sign$ by the scalar $1 - \langle
\alpha_i,s_{i-1}\dots s_1(\nu)\rangle$ and on $\bold{1}$ by the scalar $-1
-\langle\alpha_i,s_{i-1}\dots s_1(\nu)\rangle$.

\secstart\label{\Setupproof} Now fix a reducibility wall $H$ and a map $t
\mapsto \nu(t)$ as above using the normalization in \distinguishsides. We
necessarily have $R_H = \{\alpha_0, -w_0(\alpha_0)\}$ for some $\alpha_0
\in R^+$, in particular the height of $H$ is at most 2.

We apply the discussion in \Prepareproof\ to $\nu = \nu(t)$ and get
$A_{w_0}(t)$. For $t = 0$ we see that $T_i^{\alpha_i}$ (and hence $T_i$) is
invertible for all $i$ except for
$$
i \in J_H := \{j \mid s_1s_2\dots s_{j-1}\alpha_j \in R_H\}.
$$ 
Further for $i \in J_H$, $T_{i}^{\alpha_{i}}$ acts on the
$\bold{1}$-component by the scalar
$$
-\langle\alpha_i,s_{i-1}\dots s_1(\nu(t))\rangle - 1
$$
which is equal to $-2$ at $t = 0$ and it acts on
the $\sign$-component by the scalar
$$
\langle \alpha_i,s_{i-1}\dots s_1(\nu(t))\rangle - 1
$$
which is a linear function of $t$, in particular its
vanishing order is 1.

\secstart\label{\Proofofone} We keep the notation of \Prepareproof\ and
assume that we are in the situation of \Heightone, in particular $H$ is
a reducibility wall of height one and $J_H$ consists of a single element
$j$. By the arguments above we have that
$$
\ord_{t=0} \det(A_{w_0}(t)) = \dim L(\frak{M'}) = \dim E^1.
$$
Hence it follows from \orddet\ that $\omega = 1$. We will also show that in
this case we always have $\varepsilon = -1$. Let
$\lambda \in \What$ be the maximal element occuring in $L(\frak{M}')$. As the
multiplicty of $\lambda$ in $L(\frak{M}')$ is 1, the signature of the
Hermitian form induced by $\beta_1$ on the one-dimensional space $V =
\Hom_{\What}(\lambda,L(\frak{M}'))$ is equal to $\pm 1$. We have to show
that it is equal to $-1$.

By the definition we have to show that the derivative of the function
$$\rho\:t \mapsto 1 - \langle \alpha_j,s_{j-1}\dots s_1(\nu(t))\rangle$$
is negative in $t = 0$. As this is a linear function it suffices to show
that $\sgn(\rho(t)) = -\sgn(t)$. By definition of $\nu(t)$
\distinguishsides, we have $\sgn(\langle \alpha, \nu(t) \rangle - 1) =
\sgn(t)$ for the unique $\alpha \in R_H$. By \Setupproof\ we have $\alpha =
s_1s_2\dots s_{j-1}\alpha_j$. Hence we see that
$$\align
\sgn(\rho(t)) &= -\sgn(\langle \alpha_j,s_{j-1}\dots s_1(\nu(t))\rangle
-1)\\
&= -\sgn(\langle \alpha,\nu(t)\rangle)\\
&= -\sgn(t).
\endalign$$

Hence we get:

\proclaim{\sectionno Theorem}\label{\theoremone} For unramified
representations and for $\nu(0)$ lying in a reducibility wall of height
one, we have
$$\align
\lim_{t\to 0^-}\Sigma(\beta_{\nu(t)}) + \lim_{t\to 0^+}\Sigma(\beta_{\nu(t)})
&= 2\Sigma(L(\nu(0))),\\
\lim_{t\to 0^-}\Sigma(\beta_{\nu(t)}) - \lim_{t\to 0^+}\Sigma(\beta_{\nu(t)})
&= 2\Sigma(L(\frak{M}')).
\endalign$$
\endproclaim

\secstart Note that we also have
$$
\Sigma(L(\nu(0))) = \Sigma(L(\frak{M}')) \otimes \sgn
$$
as $L(\nu(0))$ is the Zelevinsky dual of $L(\frak{M}')$ \ZelonWmodule.

\secstart\label{\Proofoftwo} We keep the notation of \Setupproof\ and
assume that we are in the situation of \Heighttwo, in particular $H$ is
a reducibility wall of height two and $J_H$ consists of two elements $j_1$
and $j_2$.

If we set $m := \#W$, we have $\dim(X(\frak{M})) = m$ and it follows from
\Setupproof\ that
$$
\ord_{t=0}\det(A_{w_0}(t)) = m.
$$

If we define $\omega_1$ and $\omega_2$ as in \Heighttwo, we have by
\orddet
$$
m = \omega_1\dim(X(\frak(M)_1) + X(\frak(M)_2)) + (\omega_2 -
\omega_1)\dim(X(\frak{M}').
$$

There exist two pairs $(i_1,j_1)$ and $(i_2,j_2)$ of indices such that
$\nu_{j_k}(0) - \nu_{i_k}(0) = 1$ for $k = 1,2$. We can assume that $i_1 <
i_2$. We distinguish the cases $j_1 \not= i_2$ and $j_1 = i_2$.

By \standWmodule\ we have in the first case (resp.\ in the second case)
$$\align
&\dim(X(\frak{M}_i)) = m/2, \quad \dim(X(\frak{M}')) = m/4\\
\text{(resp.\ } &\dim(X(\frak{M}_i)) = m/2, \quad \dim(X(\frak{M}')) = m/6)
\endalign$$
and it follows that
$$\align
&m = \omega_1\frac{3}{4}m + (\omega_2 - \omega_1)\frac14 m\\
\text{(resp.\ } &m = \omega_1\frac{5}{6}m + (\omega_2 - \omega_1)\frac16 m).
\endalign$$
As $\omega_2 \geq \omega_1$ are positive integers, this implies in both cases
$$
\omega_1 = 1, \qquad \omega_2 = 2.
$$

Therefore \Jantzenfiltration\ implies
$$
\Sigma^- + \Sigma^+ = 2\Sigma(L(\nu(0))) + \varepsilon 2\Sigma(L(\frak{M}'))
$$
where
$$\align
\Sigma^+ &= \lim_{t \to 0^+} \Sigma(\beta_{\nu(t)}),\\
\Sigma^- &= \lim_{t \to 0^-} \Sigma(\beta_{\nu(t)}).
\endalign$$

Moreover it also follows that $\Sigma^- -
\Sigma^+$ is nothing but two times the signature
character of $L(\frak{M}_1) \oplus L(\frak{M}_2)$ which is zero as those
two modules are Hermitian duals of each other. 

Similarly as in \Heightone\ this implies
$$
\varepsilon = (\Sigma^+)_{\sign}\Sigma(L(\frak{M}'))_{\sgn} =
(\Sigma^-)_{\sign}\Sigma(L(\frak{M}'))_{\sgn}.
$$

Hence we get:

\proclaim{\sectionno Theorem}\label{\theoremtwo} For unramified
representations and for $\nu(0)$ lying in a reducibility wall of height
two, we have
$$\align
\Sigma^+ - \Sigma^- &= 0,\\
\Sigma^+ + \Sigma^- &= 2\Sigma(L(\nu(0))) + \varepsilon
2\Sigma(L(\frak{M}'))\\
&= 2(\Sigma(L(\frak{M}')) \otimes \sign + \varepsilon
\Sigma(L(\frak{M}')))
\endalign$$
where
$$
\varepsilon = (\Sigma^+)_{\sign}\Sigma(L(\frak{M}'))_{\sgn} =
(\Sigma^-)_{\sign}\Sigma(L(\frak{M}'))_{\sgn}.
$$
\endproclaim

\secstart Again we have
$$
\Sigma(L(\nu(0))) = \Sigma(L(\frak{M}')) \otimes \sgn
$$
as $L(\nu(0))$ is the Zelevinsky dual of $L(\frak{M}')$.


\newparagraph{Signature character for unramified representations}

\secstart\label{\Conjecture} We will now give a conjectural inductive
procedure the calculate the signature character for unramified
representations.

For this we consider Bernstein-Zelevinsky parameter giving rise to standard
an irreducible representations of $\HH^{1/2}_{GL_N}$ of the following form:
If $N$ is even, we set for $0 \leq m \leq N/2$ and real numbers $\nu_1 > \dots
> \nu_m > 0$:
$$
\frak{M}^N_{\nu_1,\dots,\nu_m} := (\nu_1,\dots,\nu_m, \undersetbrace
\text{$N/2 - m$ times} \to {[-\frac12,\frac12], \dots
,[-\frac12,\frac12]},-\nu_m,\dots,-\nu_1).
$$
If $N$ is odd, we set for $0 \leq m leq (N-1)/2$
$$
\frak{M}^N_{\nu_1,\dots,\nu_m} := (\nu_1,\dots,\nu_m,\undersetbrace
\text{$(N-1)/2-m$ times} \to {[-\frac12,\frac12], \dots,
[-\frac12,\frac12]},0,-\nu_m,\dots,-\nu_1).
$$

We conjecture the following:

\proclaim{Conjecture}: Assume that $X(\frak{M}^N_{\nu_1,\dots,\nu_m})$ is
irreducible. Then $\Sigma(X(\frak{M}^N_{\nu_1,\dots,\nu_m}))$ depends only
on the cardinality of $\{\nu_i \mid \nu_i > \frac12\}$.
\endproclaim

It follows from \theoremtwo\ that the conjecture is true in the unramified
case.

\secstart{}\label{\defineSigma} From now on we will assume that the
conjectures \Conjectureone\ and \Conjecture\ hold. In particular the
following is well-defined: For integers $N \geq 1$, $0 \leq m \leq N/2$ and
$0 \leq r \leq (N-2m)/2$ we set
$$
\Sigma^N(m,r) := \Sigma(X(\frak{M}^N_{\nu_1,\dots,\nu_m}))
$$
where $\nu_1 > \dots \nu_r > \frac12 > \nu_{r+1} > \dots > \nu_m$.

Let $\bold{e}$ be the partition $(2^m,1^{N-2m})$ of $N$. Further let $S$ be
the set of simple roots corresponding to the ordered partition
$$
\lambda := \cases (1^{N/2-m},m,1^{N/2-m}),&\text{if $N$ is even},\\
(1^{(N-1)/2-m},m+1,1^{(N-1)/2-m}),&\text{if $N$ is odd}.\endcases
$$
and let $U$ be the tempered representation of the standard Levi subgroup
corresponding to $\lambda$ which is given due to the Bernstein-Zelevinsky
classification by
$$
\undersetbrace \text{$m$ times} \to {(0) \otimes \dots \otimes (0)} \otimes
\undersetbrace \text{$N/2 - m$ times} \to {([-\frac12,\frac12], \dots
,[-\frac12,\frac12])} \otimes \undersetbrace \text{$m$ times} \to {(0)
\otimes \dots \otimes (0)}
$$
if $N$ is even and by
$$
\undersetbrace \text{$m$ times} \to {(0) \otimes \dots \otimes (0)} \otimes
\undersetbrace \text{$N/2 - m$ times} \to {([-\frac12,\frac12], \dots
,[-\frac12,\frac12])} \otimes \undersetbrace \text{$m$ times} \to {(0)
\otimes \dots \otimes (0)}
$$
if $N$ is odd.

\proclaim{Proposition} For $r = 0$ and $r = m$ we can calculate
$\Sigma^N(m,r)$ as follows:
\roster
\item"{(1)}"  For
$\bold{d} \in \hat{S}_N$ we have
$$
\Sigma^N(m,0)_{\bold d} = K_{\text{\bf d}^t,\text{\bf e}}.
$$
\item"{(2)}" Let $\Sigma_{\infty}(S,U)$ the signature at infinity
\definesigninf\ corresponding to $(S,U)$. Then we have
$$
\Sigma^N(m,m) = \Sigma_{\infty}(S,U).
$$
\endroster
\endproclaim

\demo{Proof} By \Charunitary\ we know that
$X(\frak{M}^N_{\nu_1,\dots,\nu_m})$ is unitary if all $\nu_i <
\frac12$. Hence (1) follows from \standWmodule.

By \Conjecture\ $\Sigma^N(m,m)$ is nothing but the signature at infinity
with respect to $(S,U)$. Therefore we know that the classes of
$\Sigma^N(m,m)$ and $\Sigma_{\infty}(S,U)$ in $\ZZ^{\What}/\{\pm 1\}$ are
equal. A calculation using the algorithm in \calculateinfinity\
then gives equality even in $\ZZ^{\What}$.
\qed\enddemo

\secstart\label{\inductivealgorithm} By \Conjectureone\ we have the equality
$$
\Sigma^N(m,r-1) - \Sigma^N(m,r) = 2\Sigma^N(m+1,r-1).
$$

By induction this implies for $0 \leq k \leq r$
$$
\Sigma^N(m,r) = \sum_{i=0}^k(-2)^i{k \choose i}\Sigma^N(m+i,r-k).
$$
In particular:

\proclaim{Proposition} For $N$, $m$ and $r$ as above:
$$
\Sigma^N(m,r) = \sum_{i=0}^r(-2)^i\binom ri\Sigma^N(m+i,0)
$$
which allows us to calculate $\Sigma^N(m,r)$ as we know the right hand side
by \defineSigma.
\endproclaim

\secstart\label{\charKostka} As we know not only the signature charackter
of unitary modules but also of representations at infinity, our conjectures
implies in particular: Let $w_0 \in S_N$ be the longest element and let
$\lambda$ be a partition of $N$. Set $r = [N/2]$. Then we have by
\inductivealgorithm\ and \unramifiedinfinity: 
$$\align
\chi_{\lambda}(w_0) &= \sum_{i=0}^r(-2)^i\binom ri\Sigma^N(i,0)_{\lambda}\\
&= \sum_{i=0}^r(-2)^i\binom ri K_{\lambda^t,(2^i,1^{N-2i})}.
\endalign$$


\newparagraph{Examples}

\secstart We conclude with the calculation of signature characters of
$GL_N$ for $N = 2,3,4$. Each time we will classify irreducible
$\HH_{GL_N}$-representations by the Langlands data $(S,U,\nu)$.

We will describe $S$ be the corresponding ordered partition
$(\sigma_1,\dots,\sigma_r)$ of $N$. As we consider only Hermitian modules,
we always have that $\sigma_i = \sigma_{r+1-i}$.

The tempered representation $U$ of $GL_{\sigma_1} \times \dots \times
GL_{\sigma_r}$ is a tensor product of irreducible tempered representations
$U_i$ of $GL_{\sigma_i}$ and each $U_i$ will be described by its
Bernstein-Zelevinsky datum. The condition of being Hermitian implies that
$U_i \cong U_{r+1-i}$.

Finally $\nu$ will be considered as an $r$-tuple of real numbers
$(\nu_1,\dots,\nu_r)$ with $\nu_1 > \dots > \nu_r$. We have $\nu_i +
\nu_{r+1-i} = 0$ because of the property of being Hermitian.

If $(S,U,\nu)$ is a Langlands datum, we denote by $\frak{M}(S,U,\nu)$ the
corresponding Bernstein-Zelevinsky datum.

\secstart We now consider the case $GL_2$.

$S = (2)$, $U = ([-\frac12,\frac12])$: In this case we necessarily have
$\nu = (0)$ and $\frak{M}(S,U,\nu) = ([-\frac12,\frac12])$ and $X(S,U,\nu)$
is irreducible \BZclass\ and is unitary \Charunitary. Hence by
\standWmodule:
$$
\Sigma(L(S,U,\nu)) = (1^2) = \sgn.
$$

\smallskip

$S = (2)$, $U = (0,0)$: Again $\nu$ is $(0)$, $X(S,U,\nu)$ is irreducible
and unitary, hence
$$
\Sigma(L(S,U,\nu)) = (1^2) + (2) = \sgn + {\bold 1}.
$$

\smallskip

$S = (1,1)$: $U$ is necessarily $(0) \otimes (0)$, and $\nu$ is of the form
$(\nu_1,-\nu_1)$ with $\nu_1 > 0$. We have $\frak{M}(S,U,\nu) =
(\nu_1,-\nu_1)$ and $X(S,U,\nu)$ is irreducible if and only if $\nu_1 \not=
\frac12$. For $\nu_1 < \frac12$, $L(S,U,\nu)$ is unitary, for $\nu_1 >
\frac12$, we are nearby infinity and can apply
\unramifiedinfinity. Alternatively we can use \theoremone. Finally for $\nu_1
= \frac12$ we have $L(\frac12,-\frac12) = \zeta(L([-\frac12,\frac12]))$ by
\compareZel. Hence we get
$$
\Sigma(L(S,U,(\nu_1,-\nu_1))) = \cases
(1^2) + (2),&\nu_1 < \frac12,\\
(2),&\nu_1 = \frac12,\\
-(1^2) + (2),&\nu_1 > \frac12.
\endcases
$$

\secstart We now consider the case $N = 3$:

$S = (3)$: Again $\nu = (0)$ and we have
$$
\Sigma(L(S,U,\nu)) = \cases
(1^3),&U=([-1,0,1]),\\
(1^3) + (2,1),&U = ([-\frac12,\frac12],0),\\
(1^3) + (2,1) + (3),&U = (0,0,0).
\endcases
$$

\smallskip

$S = (1,1,1)$: We have $U = (0) \otimes (0) \otimes (0)$ and $\nu =
(\nu_1,0,-\nu_1)$. The standard module $X(S,U,\nu)$ is irreducible for
$\nu_1 \not= \frac12, 1$ and it is unitary for $\nu_1 < \frac12$. At $\nu_1
= \frac12$ we have a reducibility wall of height one and at $\nu_1 = 1$ a
reducibility wall of height two. Further we have $L(\frac12,0,-\frac12) =
\zeta(L([-\frac12,\frac12],0))$ and $L(1,0,-1) = \zeta(L([-1,0,1]))$. Hence
we get
$$
\Sigma(S,U,(\nu_1,0,-\nu_1)) = \cases
(1^3) + 2(2,1) + (3),&\nu_1 < \frac12,\\
(2,1) + (3),&\nu_1 = \frac12,\\
-(1^3) +0(2,1) + (3),&\frac12 < \nu-1 < 1,\\
(3),&\nu_1 = 1,\\
(1^3) + 0(2,1) + (3),&\nu_1 > 1.
\endcases
$$

\secstart Finally consider $N = 4$:

$S = (4)$: We have $\nu = (0)$ and
$$
\Sigma(S,U,\nu) = \cases
(1^4),&U = ([-\frac32,-\frac12,\frac12,\frac32]),\\
(1^4) + (2,1^2),&U=([-1,0,1],0),\\
(1^4) + (2,1^2) + (2^2),&U=([-\frac12,\frac12],[-\frac12,\frac12]),\\
(1^4) + 2(2,1^2) + (2^2) + (3,1),&U=([-\frac12,\frac12],0,0),\\
(1^4) + 3(2,1^2) + 2(2^2) + 3(3,1) + (4),&U=(0,0,0,0).
\endcases
$$

\smallskip

$S = (2,2)$, $U = ([-\frac12,\frac12]) \otimes ([-\frac12,\frac12])$: We
have $\nu = (\nu_1,-\nu_1)$ and
$$
\frak{M}(S,U,\nu) =
([-\frac12+\nu_1,\frac12+\nu_1],[-\frac12-\nu_1,\frac12-\nu_1]).
$$
Hence $X(S,U,\nu)$ is reducible for $\nu_1 = \frac12$ and $\nu_1 = 1$, and
both are reducibility walls of height
one. For $\nu_1 < \frac12$, $L(S,U,\nu)$ is unitary and for $\nu_1 > 1$ we
are nearby infinity. Hence we can calculate both signature characters
(\standWmodule\ and \calculateinfinity ff):
$$
\Sigma(L(S,U,\nu)) = \cases
(1^4) + (2,1^2) + (2^2),&\nu_1 < \frac12,\\
(1^4) - (2,1^2) + (2^2),&\nu_1 > 1.
\endcases
$$
For $\frac12 < \nu_1 < 1$ we can use the conjecture \Conjectureone\ to
cross one of the reducibility walls. Hence the conjecture implies in this
case the following both equalities
$$\align
\Sigma(L(S,U,\nu)) &\overset \text{?} \to = (1^4) - (2,1^2) + (2^2) -
2\Sigma(L([-\frac32,-\frac12,\frac12,\frac32])),\\
\Sigma(L(S,U,\nu)) &\overset \text{?} \to = (1^4) + (2,1^2) + (2^2) -
2\Sigma(L([-1,0,1],0)).
\endalign$$
The right hand sides coincide and we get conjecturally for $\nu =
(\nu_1,-\nu_1)$ and $\frac12 < \nu_1 < 1$:
$$
\Sigma(L(S,U,\nu)) \overset \text{?} \to = -(1^4) - (2,1^2) + (2^2).
$$
Finally we can now again use \Conjectureone\ to compute
$$
\Sigma(L(S,U,(\nu_1,-\nu_1))) \overset \text{?} \to = \cases
(2^2),&\nu_1 = \frac12,\\
-(2,1^2) + (2^2),&\nu_1 = 1.
\endcases
$$

\smallskip

$S = (2,2)$, $U = (0,0) \otimes (0,0)$: We have $\frak{M}(S,U,\nu) =
(\nu_1,\nu_1,-\nu_1,-\nu_1)$ for $\nu_1 > 0$. The only reducibility wall is
at $\nu_1 = \frac12$ and it is of height 4. For $\nu_1 < \frac12$ we are
in the unitary case, for $\nu_1 > \frac12$ we are nearby infinity, and for
$\nu_1 = \frac12$ we have $L(\frac12,\frac12,-\frac12,-\frac12) =
\zeta(L([-\frac12,\frac12],[-\frac12,\frac12]))$. Hence:
$$
\Sigma(L(S,U,(\nu_1,-\nu_1))) = \cases
(1^4) + 3(2,1^2) + 2(2^2) + 3(3,1) + (4),&\nu_1 < \frac12,\\
(2^2) + (3,1) + (4),&\nu_1 = \frac12,\\
(1^4) -(2,1^2) +2(2^2) - (3,1) + (4),&\nu_1 > \frac12.
\endcases
$$

\smallskip

$S = (1,2,1)$, $U = (0) \otimes [-\frac12,\frac12] \otimes (0)$: We have
$\frak{M}(S,U,\nu) = (\nu_1,[-\frac12,\frac12],-\nu_1)$ for $\nu_1 >
0$. The reducibility walls are at $\nu_1 = \frac12$ (height one) and $\nu_1
= \frac32$ (height two). As above we get
$$
\Sigma(L(S,U,\nu)) \cases
= (1^4) + 2(2,1^2) + (2^2) + (3,1),&\nu_1 < \frac12,\\
\overset \text{?} \to = (2,1^2) + 0(2^2) + (3,1),&\nu_1 = \frac12,\\
\overset \text{?} \to = -(1^4) + 0(2,1^2) - (2^2) + (3,1),&\frac12 <
\nu_1 < \frac32,\\
\overset \text{?} \to = (2^2) - (3,1),&\nu_1 = \frac32,\\
= -(1^4) + 0(2,1^2) -(2^2) + (3,1),&\nu_1 > \frac32.
\endcases
$$

\smallskip

$S = (1,2,1)$, $U = (0) \otimes (0,0) \otimes (0)$: Here is
$\frak{M}(S,U,\nu) = (\nu_1,0,0,-\nu_1)$ for $\nu_1 > 0$. The reducibility
walls are at $\nu_1 = \frac12$ (height one) and $\nu_1 = 1$ (height 4). We
get
$$
\Sigma(L(S,U,\nu)) \cases
= (1^4) + 3(2,1^2) + 2(2^2) + 3(3,1) + (4),&\nu_1 < \frac12,\\
\overset \text{?} \to = (2,1^2) + (2,2) + 2(3,1) + (4),&\nu_1 = \frac12,\\
\overset \text{?} \to = -(1^4) - (2,1^2) + 0(2,2) + (3,1) + (4),&\frac12 <
\nu_1 < 1,\\
= (3,1) + (4),&\nu_1 = 1,\\
= -(1^4) - (2,1^2) + 0(2,2) + (3,1) + (4),&\nu_1 > 1.
\endcases
$$

\smallskip

$S = (1,1,1,1)$: Here we have $\frak{M}(S,U,\nu) =
(\nu_1,\nu_2,-\nu_2,-\nu_1)$ with $\nu_1 > \nu_2$. There are four reducibility
walls, namely those given by the conditions $\nu_1 = \frac12$ (height one),
$\nu_1 + \nu_2 = 1$ (height two), $\nu_2 = \frac12$ (height one), and
$\nu_1 - \nu_2 = 1$ (height two). We know that the signature character does
not change if we cross walls of height two \theoremtwo\ and we can use
\theoremone\ to calculate $\Sigma(L(S,U,\nu))$ for those $\nu$ such that
$X(S,U,\nu)$ is irreducible:
$$
\Sigma(L(S,U,\nu)) = \cases
(1^4) + 3(2,1^2) + 2(2^2) + 3(3,1) + (4),&\nu_2 < \nu_1 < \frac12,\\
-(1^4) - (2,1^2) + 0(2^2) + (3,1) + (4),&\nu_2 < \frac12 < \nu_1,\\
(1^4) - (2,1^2) + 2(2^2) - (3,1) + (4),&\frac12 < \nu_2 < \nu_1.
\endcases
$$
For representations lying on a single reducibility wall we
can use \theoremone\ if this wall is of height one:
$$
\Sigma(L(S,U,\nu)) = \cases
(2,1^2) + (2^2) + 2(3,1) + (4),&\nu_2 < \nu_1 = \frac12,\\
(2^2) + 0(3,1) + (4),&\frac12 = \nu_2 < \nu_1 < \frac32,\\
(2^2) + 0(3,1) + (4),&\frac12 = \nu_2 < \frac32 < \nu_1,
\endcases
$$
For representations lying on a single reducibility we use the Zelevinsky
involution if this wall is of height two:
$$
\Sigma(L(S,U,\nu)) \cases
= (2^2) + (3,1) + (4),&\nu_2 = 1 - \nu_1 < \frac12,\\
\overset \text{?} \to = -(2^2) + (3,1) + (4),&\nu_2 = \nu_1 - 1 < \frac12,\\
= (2^2) - (3,1) + (4),&\nu_2 = \nu_1 - 1 > \frac12.
\endcases
$$
Finally there is a unique representation which lies on two reducibility
walls, namely $L(\frac32,\frac12,-\frac12,-\frac32)$ and this
representation is unitary. It is the Zelevinsky dual of
$L([-\frac32,-\frac12,\frac12,\frac32])$ and hence we have in this case
$$
\Sigma(L(S,U,\nu)) = (4).
$$



\vskip 1cm \Refs\nofrills{References}


\vskip 1cm


\widestnumber\key{LuCUSPIII}


\ref \key BM1 \by D. Barbasch and A. Moy  \paper A unitarity criterion
for $p$-adic groups \jour Invet. math. \vol 98 \pages  19--37 \yr
1989\endref


\ref \key BM2 \by   D. Barbasch and A. Moy \paper Reduction to real
infinitesimal character in affine Hecke algebras\jour Jour. of  the
AMS\vol 6 \#3 \pages  611--635 \yr 1993\endref


\ref \key BM3 \by  D. Barbasch and A. Moy \paper Unitary Spherical
Spectrum for $p$-adic classical groups \jour Acta Applicandae
Math.\vol 44 \pages 3--37 \yr 1996\endref


\ref \key Ca \by P.~Cartier \paper Representations of $p$-adic groups:
a survey \jour Proc. of Symp. in Pure Math. \vol 33 (1) \pages
111--155 \yr 1979 \endref


\ref \key CP \by C.~de~Concini, C.~Procesi \paper Symmetric Functions,
Conjugacy Classes and the Flag Variety \jour Inv. Math. \vol 64 \pages
203--219 \yr 1981 \endref


\ref \key Ev \by S. Evens \paper The Langlands classification for
graded Hecke algebras \jour Proc. of the AMS \vol 124 \pages
1285--1290 \yr 1996 \endref


\ref \key KL\by D. Khazhdan and  G. Lusztig   \paper Proof of
Deligne-Langlands conjecture for Hecke algebras \jour Inv. math.\vol
87 \pages 153--215 \yr 1987\endref


\ref \key Ku\by S. Kudla \paper The Local Langlands Correspondence:
The Non-Archimedian case \jour Proc. of Symp. Pure. Math. \vol 55 (2)
\pages 365--391 \year 1994 \endref


\ref \key KW\by  A. Kent and G. Watts \paper Signature Character for
$A_2$ and $B_2$ \jour Commun. Math. Phys.\vol  143\pages 1--16 \yr
1991 \endref


\ref \key KZ\by H. Knight, A. Zelevinsky \paper Representations of
quivers of type $A$ and the multisegment duality \jour Adv. Math. \vol
117 \yr 1996 \pages 273--293 \endref


\ref\key Lu1 \by G.~Lusztig \paper Cuspidal Local systems and
Graded Hecke Algebras I \jour Publ. Math. Institut des Hautes \'Etudes
Scientifiques \vol 67 \yr 1988 \pages 145--202 \endref


\ref \key Lu2\by G. Lusztig  \paper Affine Hecke algebras and their
graded version \jour J. AMS\vol 2 \pages 599--635 \yr 1989\endref
 
\ref\key Lu3 \by G. Lusztig   \paper Classification of Unipotent
Representations of simple $p$-adic groups \jour IMNR\vol 11 \pages
517-589 \yr 1995\endref


\ref\key Lu4 \by G. Lusztig \paper Cuspidal Local systems and
Graded Hecke Algebras II \inbook Representations of groups
(ed. B.~Allison and G.~Cliff) Canad. Math. Soc. Conf. Proc. \vol 16
\yr 1995 \publ Amer. Math. Soc. \pages 217--275 \endref


\ref\key Lu5 \by G.~Lusztig \paper Cuspidal Local systems and
Graded Hecke Algebras III \paperinfo preprint MIT, RT/0108173 \yr 2001
\endref


\ref\key Ma \by I.G.~Macdonald \book Symmetric Functions and Hall
Polynomials \bookinfo 2nd edition \publ Oxford University press \yr
1995 \endref

 
\ref \key MW\by C.~Moeglin, J.-L.~Waldspurger\paper Sur l'involution de
Zelevinski \jour J. Reine Angew. Math. \vol 372 \yr 1986 \pages 136--177\endref


\ref \key Ta \by M. Tadi\v c  \paper Classification of unitary
representations in irreducible representation of general linear group
(non-archimedian case) \jour Ann. scient. Ec. Norm. Sup.\vol 19 \pages
335-382 \yr 1986 \endref


\ref \key Vo \by D.~A.~Vogan, Jr. \paper Unitarizability of certain series
of representations \jour Annals of Math.\vol 120 \pages 141-187 \yr 1984
\endref


\ref \key Ze \by A.V. Zelevinsky \paper Induced representations of
reductive $p$-adic groups II, on irreducible representations of
$GL(n)$ \jour Ann. Sci. ENS $4^e$ s\'erie \vol 13 \pages 165--210 \yr
1980 \endref




 


\endRefs
\enddocument